\documentclass[10pt, a4paper]{article}
\usepackage{amsmath}
\usepackage{amssymb}
\usepackage{graphicx}
\usepackage{psfrag}
\usepackage{latexsym, epsfig}

\newtheorem{theorem}{Theorem}[section]
\newtheorem{lemma}[theorem]{Lemma}%[chapter]
\newtheorem{proposition}[theorem]{Proposition}%[chapter]
\newtheorem{definition}[theorem]{Definition}%[chapter]
\newtheorem{remark}[theorem]{Remark}%[chapter]
\newcommand{\conv}{\hbox{\rm conv}}

\newcommand{\conc}{\hbox{\rm conc}}

\newcommand{\elleuno}{\mathbb{L}^1}

\newcommand{\real}{\mathbb{R}}
\newcommand{\nat}{\mathbb{N}}
\newcommand{\M}{\mathbf{M}}

\newcommand{\Q}{\mathcal{Q}}
\newcommand{\RC}{\mathcal R}
\newcommand{\SC}{\mathcal S}
\newcommand{\MC}{\mathcal{M}}
\newcommand{\U}{\mathcal{U}}
\newcommand{\integer}{\mathbb{Z}}
\newcommand{\OL}{\mathcal{O}}
\newcommand{\C}{\mathcal{C}}

\newcommand{\J}{\mathcal{J}}

\newcommand{\eps}{\varepsilon}
\newcommand{\vfi}{\varphi}
\newcommand{\ups}{\Upsilon}
\newcommand{\TV}{\hbox{\rm Tot.Var.}}

\newcommand{\dtau}{\displaystyle{\frac{d}{d\tau}}}
\newcommand{\second}{{\prime\prime}}

\newcommand{\scalar}[2]{{\big\langle #1, #2 \big\rangle}}

\newcommand{\wl}{\widetilde{\lambda}}
\newcommand{\wwl}{\widetilde{\widetilde{\lambda}}}

\newcommand{\wy}{\widetilde{y}}
\newcommand{\wwy}{\widetilde{\widetilde{y}}}
\newcommand{\wwyhk}{\widetilde{\widetilde{y}}^h_k}
\newcommand{\ws}{\widetilde{s}}
\newcommand{\wws}{\widetilde{\widetilde{s}}}

\newcommand{\qed}{{\hfill $\square$}\medskip}

\begin{document}
\title{A locally quadratic Glimm functional and sharp convergence
  rate of the Glimm scheme for nonlinear hyperbolic systems}
\author{Fabio Ancona\footnote{Dipartimento di Matematica and C.I.R.A.M..
Via Saragozza 8, 40123 - Bologna, Italy,
email: ancona@ciram.unibo.it and ancona@math.unipd.it}
\and Andrea Marson\footnote{Dipartimento di Matematica Pura ed Applicata,
Via Trieste 63, 35121 - Padova, Italy,
email: marson@math.unipd.it}}
\maketitle
\begin{center}
\date{September 1\ 2008}
\end{center}
\begin{abstract}
Consider the Cauchy problem for a strictly hyperbolic,  $N\times N$
 quasilinear system in one space dimension
$$
u_t+A(u) u_x=0,\qquad u(0,x)=\bar u(x),
\eqno (1)
$$
where $u \mapsto A(u)$
is a smooth matrix-valued map, and the initial data $\overline u$ is
assumed to have small total variation. We investigate the rate of
convergence of approximate solutions of (1) constructed by
the Glimm scheme, under the assumption that, 
letting $\lambda_k(u)$, $r_k(u)$ denote the $k$-th eigenvalue and a
corresponding eigenvector of $A(u)$, respectively, for each $k$-th
characteristic family the linearly degenerate manifold
$$
 \mathcal{M}_k \doteq \big\{u\in\Omega~:~\nabla\lambda_k(u)\cdot
 r_k(u)=0\big\}
$$
is either the whole space, or
it is empty, or it
consists of a finite number of smooth, $N\!-\!1$-dimensional,
connected, manifolds that are transversal to the characteristic vector
field $r_k$\,.
We introduce a Glimm type functional which is the sum of the 
cubic interaction potential defined in \cite{sie}, and of a quadratic term
that takes into account 
interactions of waves of the same family with strength smaller than some
fixed threshold parameter.
Relying on an adapted wave tracing method,
and on the decrease amount
of such a functional, we obtain the same
type of error estimates valid for Glimm approximate solutions of
hyperbolic systems satisfying the classical Lax assumptions of genuine
nonlinea\-ri\-ty or linear degeneracy of the characteristic families.
\end{abstract}
\newpage
             
%%%%%%%%%%%%%%%%%%%%%%%%%%%%%%%%%%%%%%%%%%%%%%%%%%%%%%%%%%%%%%%%%%%%%%
%
%
%	INTRODUCTION
%
%

\tableofcontents

\section[Introduction]{Introduction}
\label{sec:int}

Consider the Cauchy problem for a general system of hyperbolic
conservation laws in one space dimension
\begin{align}
  \label{syscon}
  \null
  & u_t+F(u)_x =0\,,\\
  \noalign{\smallskip}
  \label{inda}
  & u(0,x) =\overline u(x)\,.
\end{align}
Here the vector $u=u(t,x)=\big(u_1(t,x),\dots, u_N(t,x)\big)$
represents the {\it conserved quantities}, while the components of the
vector valued function
$$
F(u)=\big(F_1(u),\dots,F_N(u)\big)
$$
are the
corresponding {\it fluxes}.  We assume that the flux function $F$ is a
smooth map defined on a domain $\Omega\subseteq\real^N$, and that the
system (\ref{syscon}) is strictly hyperbolic, i.e. that the Jacobian
matrix $A(u)=DF(u)$ has $N$ real distinct eigenvalues
\begin{equation}
  \lambda_1(u)<\dots<\lambda_N(u)\qquad\forall~u\,.
  \label{strhyp}
\end{equation}
Denote with $r_1(u),\dots,r_N(u)$ a corresponding basis of right
eigenvectors.  Hyperbolic equations in conservation form physically
arise in several contexts. A primary example of such systems is
provided by the Euler equations of non-viscous gases, see \cite{Daf}.

It is well known that, because of the nonlinear dependence of the
characteristic speeds $\lambda_k(u)$ on the state variable $u$,
classical solutions to (\ref{syscon}) can develop discontinuities
(shock wave) in finite time, no matter of the regularity of the
initial data.  Therefore, in order to construct solutions globally
defined in time, one must consider weak solutions interpreting the
equation (\ref{syscon}) in a distributional sense.  Moreover, for sake
of uniqueness, an entropy criterion for admissibility is usually added
to rule out non-physical discontinuities.  In \cite{tplrp2p2,tplrpnpn}
T.P.~Liu proposed the following admissibility criterion valid for weak
solutions to general systems of conservation laws, that generalizes
the classical stability condition introduced by Lax \cite{lax}.
\begin{definition}
\label{def:Laxstab}
A shock discontinuity of the $k$-th family $(u^L,\,u^R)$, traveling
with speed $\sigma_k[u^L,\,u^R]$, is {\it Liu admissible} if,
for any state $u$ lying on the Hugoniot curve $S_k[u^L]$ between
$u^L$ and $u^R$,
the shock speed $\sigma_k[u^L,u]$ of the discontinuity $(u^L,u)$ satisfies
\begin{equation}
  \label{liuc}
  \sigma_k[u^L,\,u]\geq\sigma_k[u^L,\,u^R]\,.
\end{equation}
\end{definition}

The existence of global weak admissible solutions to
(\ref{syscon})-(\ref{inda}) with small total variation was first
established in the celebrated paper of Glimm~\cite{glimm}
under the additional assumption that each characteristic field~$r_k$
be either {\it linearly degenerate} (LD), so that
\begin{equation}
  \nabla\lambda_k(u)\cdot r_k(u)=0\qquad\forall~u\,,
\end{equation}
or else {\it genuinely nonlinear} (GNL) i.e.
\begin{equation}
  \nabla\lambda_k(u)\cdot r_k(u)\neq 0\qquad\forall~u\,.
\end{equation}
A random choice method, the Glimm scheme, was introduced in
\cite{glimm} to construct approximate solutions of the general Cauchy
problem (\ref{syscon})-(\ref{inda}) by piecing together solutions of
several {\it Riemann problems}, i.e. Cauchy problems whose initial
data are piecewise constant with a single jump at the origin
\begin{equation}
  u(0,x)=\begin{cases} u^L\qquad &\text{if\quad $x<0$\,,}\\
  u^R\qquad &\text{if\quad $x>0$\,.}
  \end{cases}
  \label{riemin}
\end{equation}
Using a nonlinear functional introduced by Glimm, that measures the
nonlinear coupling of waves in the solution, one can establish
a-priori bounds on the total variation of a family of approximate
solutions.
These uniform estimates then yield the convergence of a sequence of
approximate solutions to the weak admissible solution of
(\ref{syscon})-(\ref{inda}).  The existence theory for the Cauchy
problem (\ref{syscon})-(\ref{inda}) based on a Glimm scheme was
extended by Liu~\cite{tplams}, Liu and Yang~\cite{glimmly}, and by
Iguchy and LeFloch~\cite{IL} to the case of systems 
with {\it non genuinely nonlinear} (NGNL) characteristic
families  whose flux
function satisfy the more general assumption:
\begin{description}
\item[{\rm (H)}] 
  The vector valued function $F$ is $\C^3$, and
  for each $k\in\{1,\dots,N\}$-th characteristic family
  the linearly degenerate manifold
  \begin{equation}
    \mathcal{M}_k \doteq \big\{u\in\Omega~:~\nabla\lambda_k(u)\cdot
    r_k(u)=0\big\}
    \label{LDM}
  \end{equation}
  is either empty (GNL characteristic field),
  or it is the whole space (LD characteristic field), or it consists of a finite number $\leq M$ of smooth,
  $N\!-\!1$-dimensional, connected, manifolds, and there holds
\begin{equation}
  \label{eq:transv-ld}
   \nabla (\nabla\lambda_k \cdot r_k)(u)\cdot r_k(u)\neq 0\qquad\forall u\in \mathcal{M}_k\,.
\end{equation}
\end{description}
\vspace{.2truecm}

Aim of the present paper is to provide a sharp convergence rate for
approximate solutions obtained by the Glimm scheme valid for strictly
hyperbolic systems of conservation laws satisfying the assumption (H).
We recall that in the Glimm scheme, one works with a fix grid in the
$t$-$x$ plane, with mesh sizes $\Delta t, \Delta x$.  An approximate
solution $u^\eps$ of \eqref{syscon}-(\ref{inda}) is then constructed
as follows.  By possibly performing a linear change of coordinates in
the $t$-$x$ plane, we may assume that the characteristic speeds
$\lambda_k(u)$, $1\leq k\leq N$, take values in the interval $[0,1]$,
for all $u\in\Omega$.  Then, choose $\Delta t=\Delta x\doteq \eps$,
and let $\{ \theta_\ell \}_{\ell\in \nat}\subset [0,1]$ be an
equidistributed sequence of numbers, which thus satisfies the
condition
\begin{equation}
 \label{eq:eqdistr}
\lim_{n\to\infty}
\bigg|
\lambda-\frac{1}{n}\sum_{\ell=0}^{n-1} \chi_{[0,\lambda]} (\theta_\ell)
\bigg|
=0\qquad\forall~\lambda\in[0,1]\,,
\end{equation}
where $\chi_{[0,\lambda]}$ denotes the characteristic function of
the interval $[0,\lambda]$. On the initial strip $0\leq t<\eps$,
$u^\eps$ 
is defined as the exact solution of \eqref{syscon}, with
starting condition
\begin{equation}
\nonumber
u^\eps(0,x)=\overline u\big((j+\theta_0)\eps\big)\qquad
\forall~x\in\,]j\eps,\,(j+1)\eps\,[\,.
\end{equation}
Next, assuming that $u^\eps$ has been constructed for $t\in
[0,\,i\eps[$\,, on the strip $i\eps\leq t <(i+1)\eps$\,, $u^\eps$ is
defined as the exact solution of \eqref{syscon}, with starting
condition
\begin{equation}
\nonumber
u^\eps(i\eps,x)=u^\eps\big(i\eps-,\,(j+\theta_i)\eps\big)\qquad
\forall~x\in\,]j\eps,\,(j+1)\eps\,[\,.
\end{equation}
Relying on uniform a-priori bounds on the total variation, we thus
define inductively the approximate solution $u^\eps(t,\cdot)$ for all
$t\geq 0$.

One can repeat this construction with the same values $\theta_i$ for
each time interval $[i\eps,\,(i+1)\eps[\,$, and letting the mesh size
$\eps$ tend to zero. Hence, we obtain a sequence of approximate
solutions
which converge, by compactness, to some limit function $u$ that is
shown to be a weak admissible solution of~\eqref{syscon}-(\ref{inda})
(cfr.~\cite{tplGlimm}). In order to derive an accurate estimate of the
convergence rate of the approximate solutions, it was introduced
in~\cite{bm} an equidistributed sequence $\{ \theta_\ell \}_{\ell\in
\nat}\subset [0,1]$ enjoying the following property. For any given
$0\leq m<n$, define the \textit{discrepancy} of the set $\big\{
\theta_m,\dots, \theta_{n-1} \big\}$ as
\begin{equation}
  \label{eq:dis}
  D_{m,n} \doteq \sup_{\lambda\in [0,1]} \left\vert \lambda -
  \frac{1}{n-m} \sum_{m\leq\ell <n} \chi_{[0,\lambda]} (\theta_\ell)
  \right\vert \,.
\end{equation}
Then, there holds
\begin{equation}
\label{eq:disest}
D_{m,n} \leq \mathcal{O}(1) \cdot \frac{1+\log(n-m)}{n-m} \qquad
\forall~n>m\geq 1\,.
\end{equation}
Here, and throughout the paper, $\mathcal{O}(1)$ denotes a uniformly
bounded quantity, while we will use the Landau symbol $o(1)$ to
indicate a quantity that approaches zero as $\eps\to 0$.  Relying on
the existence of a Lipschitz continuous semigroup of solutions
generated by (\ref{syscon}), compatible with the solutions of the
Riemann problems,
it was proved in~\cite{bm} that, for systems with GNL or LD
characteristic fields, the $\elleuno$ convergence rate of the Glimm
approximate solutions constructed in connection with a sequence
enjoying the property~\eqref{eq:disest} is $\text{o}(1)\cdot
\sqrt\eps\,|\ln\eps|$.  In the case of general systems satisfying the
assumption (H), it was derived in~\cite{HY} an estimate of the
$\elleuno$ norm of the error in the Glimm approximate solutions of the
order $\text{o}(1)\sqrt[3]\eps\,|\ln\eps|$.

In the present paper, we improve this result by establishing the same
convergence rate of the approximate solutions generated by the Glimm
scheme for systems satisfying the assumption (H) as in the case of
systems with GNL or LD characteristic fields.  Namely, our result is
the following.
\begin{theorem}
\label{thm:glimm}
Let $F$ be a $\C^3$ map from a domain $\Omega\subseteq\real^N$ into
$\real^N$ satisfying the assumption (H), and assume that the system (\ref{syscon}) is strictly
hyperbolic. Given an initial datum $\overline{u}$ with small total variation,
let $u(t,\cdot)$
be the unique Liu admissible solution of
(\ref{syscon})-(\ref{inda}).
Let $\{ \theta_k \}_{k\in \nat}\subset [0,1]$ be a sequence satisfying
(\ref{eq:disest}) and construct the corresponding Glimm approximate
solution $u^\eps$ of (\ref{syscon})-(\ref{inda}) with mesh sizes
$\Delta x=\Delta t=\eps$. Then, for every $T\geq 0$ there holds
\begin{equation}
\label{eq:Glimm}
\lim_{\eps\to 0} \frac{\Vert u^\eps (T,\cdot) - u(T,\cdot) \Vert_{\elleuno}}
{\sqrt{\eps} \vert \log \eps \vert}\,,
\end{equation}
and the limit is uniform w.r.t. $\overline{u}$ as long as
$\TV (\overline{u})$ remains uniformly small.
\end{theorem}

Our result applies more generally to strictly hyperbolic $N\times N$
quasilinear systems
\begin{equation}
  u_t+A(u)\, u_x =0\,,
  \label{sys}
\end{equation}
not necessarily in conservation form, where $A$ is a $\C^2$ matrix
valued map defined from a domain $\Omega\subseteq\real^N$ into
$\M^{N\times N}(\real)$, whose eigenvalues $\lambda_k,$
$k\in\{1,\dots,N\}$, satisfy the assumption stated in (H).
Indeed, one may alternatively assume that 
$A:\Omega\to \M^{N\times N}(\real)$ is a $\C^{1,1}$   map,
and that for each NGNL $k\in\{1,\dots,N\}$-th characteristic family
  the linearly degenerate manifold $\mathcal{M}_k$ consists of
a finite number of 
connected manifolds $\mathcal{M}_{k,h}$, 
that are either $N\!-\!1$-dimensional as in (H)
or $N$-dimensional
with a similar 
condition to \eqref{eq:transv-ld} 
(cfr. Remark~\ref{rem3:c11ngnl} in \S~\ref{sec:oip} and Remark~\ref{rem8:genngnl}
in \S~\ref{sec:con}).

In fact, the
fundamental paper of Bianchini and Bressan \cite{BB} shows that,
for any $\C^{1,1}$  map $A:\Omega\to \M^{N\times N}(\real)$ with strictly hyperbolic values,
(\ref{sys}) generates a unique (up to the domain) Lipschitz continuous
semigroup $\{S_t\,:\ t\geq 0\}$ of {\it vanishing viscosity solutions}
obtained as the (unique) limits of solutions to the (artificial)
viscous parabolic approximation
\begin{equation}
 \label{eq1:vanvisnc}
    u_t+A(u)\, u_x = \mu\, u_{xx}\,,
    \eqno{(1.16)_\mu}
    \nonumber\setcounter{equation}{16}
\end{equation}
when the viscosity coefficient $\mu\to 0$.
The trajectories of such a semigroup starting
from piecewise constant initial data locally coincide with the {\it
``admissible'' solution} of each Riemann problem determined by the
jumps in the initial data. Moreover,
any limit of Glimm approximations coincides with the
corresponding trajectory of the semigroup generated by \eqref{sys}.
In particular, in the conservative case
where $A(u)=DF(u)$ every vanishing viscosity solution of the Cauchy
problem~(\ref{sys})-(\ref{inda}) provides a weak solution
of~(\ref{syscon})-(\ref{inda}) satisfying the Liu admissibility
conditions~\eqref{liuc}. 
\vspace{.2truecm}

The proof of the error bound \eqref{eq:Glimm} follows the same
strategy adopted in~\cite{bm}, relying on the careful analysis of the
structure of the solution for systems satisfying the assumption (H),
developed by T.P.~Liu and T.~Yang in~\cite{tplams,glimmly}.  Indeed,
to estimate the distance between a Lipschitz continuous (in time)
approximate solutions $w$ of \eqref{sys} and the corresponding exact
solution one would like to use the error bound~\cite{lnb}
\begin{equation}
\label{smgrbound}
\big\|w(T)-S_tw(0)\big\|_{\elleuno}\leq L \int_0^T
\liminf_{h\to 0+}\frac{\big\|w(t+h)-S_h w(t)\big\|_{\elleuno}}{h}~dt\,,
\end{equation}
where $L$ denotes a Lipschitz constant of the semigroup $S$ generated
by~\eqref{syscon}.  However, for approximate solutions constructed by
the Glimm scheme, a direct application of this formula is of little
help because of the additional errors introduced by the restarting
procedures at times $t_i\doteq i \eps$.  For this reason, following
the wave tracing analysis in~\cite{glimmly}, it is useful to partition
the elementary waves present in the approximate solution, say in a
time interval $[\tau_1, \tau_2]$, into virtual waves that can be
either traced back from $\tau_2$ to $\tau_1$, or are canceled or
generated by interactions occurring in $[\tau_1, \tau_2]$.  Thanks to
the simplified wave pattern associated to this partition, one can
construct a front tracking approximation having the same initial and
terminal values as the Glimm approximation, and thus
establish~\eqref{eq:Glimm} relying on~\eqref{smgrbound}.

As one would expect, the presence of elementary waves with various
composite wave patterns for systems satisfying the assumption (H),
requires a careful analysis of the errors introduced by this
wave-partition algorithm.  As customary, the change of wave-size and
wave-speeds when an interaction takes place is controlled by a Glimm
functional that measures the potential interaction of waves in the
solution.

For general strictly hyperbolic systems \eqref{sys} satisfying the
assumption (H), several nonlinear functionals were introduced
in~\cite{tplams,glimmly,IL,sie}, consisting of a standard Glimm
quadratic functional, for the interaction of waves of different
families, and of a cubic functional measuring the potential
interaction between waves of the same family. This cubic part of the
functional is defined in terms of the strengths of any pair of waves
of the same family and of the absolute value of the angle between
them~\cite{sie} (or of the positive part of the angle between two
waves~\cite{tplams,glimmly}). Such functionals work perfectly to
establish uniform a-priori bounds on the total variation of the
solution, but are not effective to control the quadratic order error
produced by the change of wave speeds for interactions of waves of the
same family, of arbitrarily small sizes.

On the other hand, in the case of systems whose characteristic
families admit a single, connected, $N-1$-dimensional degenerate
manifold \eqref{LDM}, it was introduced in \cite{amwftnpn} a
decreasing potential interaction functional which is of second order
w.r.t. the total variation (measuring the potential interaction
between {\it any} pair of waves as proportional to the product of
their strengths, no matter if they belong to the same family or not).

In the present paper, in connection with a fixed threshold 
parameter $\delta_0>0$, we define a Glimm type functional
$Q\doteq Q_q+ c\, \Q$, for a suitable constant $c>0$,
which is the sum of a quadratic term $Q_q$ and of the cubic
interaction potential~$\Q$ defined in \cite{sie}.  Here, in presence
of interactions between waves of the same families and strength
smaller than $\delta_0$, $Q_q$ behaves as the
interaction functional introduced in \cite{amwftnpn}, while the
decrease of~$\Q$ controls the possible increase of $Q_q$ at
interactions involving waves of the same family and strength larger
than $\delta_0$.  Employing this functional we can produce a
simplified wave partition pattern whose errors are controlled by the
total decrease of the Glimm functional in the time interval taken in
consideration, and thus yield the error estimate \eqref{eq:Glimm}.
\smallskip

{\bf Note added.} During the completion of the present paper, 
we have had knowledge of a contemporary different proof of
the same convergence rate~\eqref{eq:Glimm} provided by
J.~Hua, Z.~Jiang and T.~Yang~\cite{HJY}, for Glimm approximations of
a system~\eqref{sys} satisfying the assumption (H). Their proof is
obtained by using an adapted form of the functional introduced
in~\cite{tplams,glimmly}, that takes care of the errors in the
wave-speeds at interactions between waves of the same family. 

%
%
%	PRELIMINARIES
%
%
\section[Preliminaries]{Preliminaries}
\label{sec:pre}
\setcounter{equation}{0}

Let $A$ be a smooth matrix-valued map defined on a domain
$\Omega\subset\real^N$, with values in the set of $N\times N$
matrices.  Assume that each $A(u)$ is strictly hyperbolic and denote
by $\{\lambda_1(u),\,\dots,\,\lambda_N(u)\}\subset [0,1]$ its
eigenvalues.  Since we will consider only solutions with small total
variation that take values in a neighborhood of a compact set
$K\subset\Omega$, it is not restrictive to assume that $\Omega$ is
bounded and that there exist constants $\widehat\lambda_0<\cdots<
\widehat\lambda_{N}$ such that
\begin{equation}
  \label{hyp}
  \widehat\lambda_{k-1}<\lambda_k(u)<\widehat\lambda_k\,,\qquad
  \forall~u\,,\quad k=1,\dots,N\,.
\end{equation}
One can choose bases of right and left eigenvectors $r_k(u)$,
$l_k(u)$, $(k=1,\dots,N)$, associated to $\lambda_k(u)$, normalized so
that
\begin{equation}
  \label{eq2:norm}
  \big\vert r_k(u) \big\vert \equiv 1\,, \qquad \big\langle
  l_h(u),\,r_k(u) \big\rangle = 
  \begin{cases}
    1\quad &\text{if\quad $k=h\,,$}\\
    0\quad &\text{if\quad $k\neq h\,,$}
  \end{cases}
  \qquad\forall~u\,.
\end{equation}
By the strict hyperbolicity of the system, in the conservative case~\eqref{syscon}
(where $A(u)=DF(u)$), for every fixed
$u_0\in\Omega$ and for each $k\in\{1,\dots,N\}$-th characteristic
family one can construct in a neighborhood of $u_0$ a one-parameter
smooth curve $S_k[u_0]$ passing through $u_0$ (called the $k$-th {\it
Hugoniot curve} issuing from $u_0$), whose points $u\in S_k[u_0]$
satisfy the Rankine Hugoniot equation $F(u)-F(u_0)=\sigma
\big(u-u_0\big)$\, for some scalar~$\sigma=\sigma_k[u_0,u]$.  The
curve $S_k[u_0]$ is tangent at~$u_0$ to the right eigenvector
$r_k(u_0)$ of~$A(u_0)$ associated to $\lambda_k(u_0)$, and we say that
$(u^L,\,u^R)$ is a {\it shock discontinuity} of the $k$-th family with
{\it speed} $\sigma_k[u^L,\,u^R]$ if~$u^R\in S_k[u^L]$.
\vspace{.2truecm}

We describe here the general method introduced in \cite{BB,srp} to
construct the self-similar solution of a Riemann problem for a
strictly hyperbolic quasilinear system (\ref{sys}).  As customary,
the basic step consists in constructing the {\it elementary curve of
the $k$-th family} $(k=1,\dots,N)$ for every given left state $u^L$,
which is a one parameter curve of right states $s\mapsto T_k[u^L](s)$
with the property that the Riemann problem having initial data $(u^L,
u^R)$, $u^R\doteq T_k[u^L](s)$, admits a vanishing viscosity solution
consisting only of {\it elementary waves} of the $k$-th characteristic
family.  Such a curve is constructed by looking at the fixed point of
a suitable contractive transformation associated to a smooth manifold
of viscous traveling profiles for the parabolic system with unit
viscosity $(1.16)_1$.

Given a fixed state $u_0\in\Omega$, and an index $k\in\{1,\dots,N\},$
in connection with the $N+2$-dimensional smooth manifold of bounded 
traveling profiles of $(1.16)_1$ with speed close to $\lambda_k(u_0)$,
one can define on a neighborhood of
$(u_0,0,\lambda_k(u_0)) \in \real^N\times\real\times\real$
suitable smooth vector functions $(u,v_k,\sigma) \mapsto\widetilde
r_k (u,v_k,\sigma)$ that
satisfy
\begin{equation}
  \label{eq2:ri}
  \widetilde r_k \big(u_0,0,\sigma\big) = r_k(u_0) \qquad
  \forall~\sigma\,,
\end{equation}
and are normalized so that
\begin{equation}
  \label{normrt}
  \scalar{l_k(u_0)}{\widetilde r_k (u,v_k,\sigma)} = 1 \qquad
  \forall~u\,,\,v_k\,,\,\sigma\,.
\end{equation}
The vector valued map
$\widetilde r_k(u,v_k,\sigma)$ is called the $k$-th \emph{generalized
eigenvector} of the matrix~$A(u)$, associated to the
\emph{generalized eigenvalue}
$$\widetilde \lambda_k(u,v_k,\sigma)\doteq\scalar{l_k(u_0)}{A(u)\,
    \widetilde r_k (u,v_k,\sigma)}\,,
$$
that satisfies  the identity
\begin{equation}
\label{eq2:lk}
  \widetilde \lambda_k \big(u_0,v_k,\sigma \big) =
  \lambda_k(u_0)\qquad 
  \forall~v_k\,,\,\sigma\,.
\end{equation}
Next, given a left state $u^L$ in a neighborhood of $u_0$ and
$0<s<<1$, consider 
the integral system
\begin{equation}
  \label{eq2:Ti}
  \begin{cases}
    u(\tau) = u^L + \displaystyle{\int_0^\tau} \widetilde r_k
    \big( u(\xi), v_k(\xi), \sigma(\xi)
    \big)~d\xi\,,\\
    \noalign{\medskip}
    v_k(\tau) = \widetilde F_k\big(\tau;\,u,v_k,\sigma\big) -
    \conv_{[0,s]} \widetilde F_k\big(\tau;\,u,v_k,\sigma\big)\,,\\
    \noalign{\medskip}
    \sigma(\tau) = \dtau \conv_{[0,s]} \widetilde
    F_k\big(\tau;\,u,v_k,\sigma\big)\,,
  \end{cases}
  \qquad
  0\leq\tau \leq s\,,
\end{equation}
where $\tau\mapsto \widetilde f_k(\tau)\doteq\widetilde F_k(\tau;\,u,v_k,\sigma)$ is the
{\it``reduced flux function''} associated to (1.16) defined, by
\begin{equation}
  \label{eq2:tFi}
   \widetilde f_k(\tau)\doteq \int_0^\tau \widetilde
  \lambda_k \big( u(\xi), v_k(\xi), \sigma(\xi) \big)~d\xi\,,
\end{equation}
and we let $\conv_{[0,s]} \widetilde f_k(\tau)$ denote the
{\it lower convex envelope} of~$\widetilde f_k$ on~$[0,s]$,  i.e.
\begin{multline}
\label{convenv}
  \text{conv}_{[0,s]}\widetilde f_k(\tau)\doteq
  \inf\Big\{
  \theta\,\widetilde f_k(y)+(1-\theta)\,\widetilde f_k(z)\,:\\
  \theta\in[0,1]\,, \
  y,\,z\in[0,s]\,, \ \tau=\theta y+(1-\theta)z
  \Big\}\,.
\end{multline}
Relying on (\ref{eq2:ri}), (\ref{eq2:lk}) it is shown
in~\cite{BB,srp} that, for $s$ sufficiently small, the transformation
defined by the right-hand side of~(\ref{eq2:Ti}) maps a domain of
continuous curves $\tau\mapsto (u(\tau),v_k(\tau),\sigma(\tau))$ into
itself, and is a contraction w.r.t. a suitable weighted norm.  Hence,
for every $u^L$ in a neighborhood $\U_0$ of $u_0$, the transformation
defined by~(\ref{eq2:Ti}) admits a unique fixed point
\begin{equation}
  \tau\mapsto
  \big(u(\tau;\,u^L,s),\ v_k(\tau;\,u^L,s),\ \sigma(\tau;\,u^L,s)\big)
  \qquad\tau\in[0,s]\,,
  \nonumber
\end{equation}
which provides a Lipschitz continuous solution to the integral system
(\ref{eq2:Ti}).  The elementary curve of right states of the $k$-th
family issuing from $u^L$ is then defined as the terminal value at
$\tau=s$ of the $u$-component of the solution to the integral system
(\ref{eq2:Ti}), i.e. by setting
\begin{equation}
  \label{Tdef}
  T_k[u^L](s)\doteq u(s;\,u^L,s)\,.
\end{equation}
Sometimes, the value \eqref{Tdef} of the elementary curve issuing from $u^L$
will be equivalently 
written $T_k(s)[u^L]$.
In the following it will be convenient to adopt the notations
\begin{equation}
  \arraycolsep 2pt
  \begin{array}{rl}
    \sigma_k[u^L](s,\tau) &\doteq \sigma(\tau;\,u^L,s)
    \\
    \noalign{\medskip}
    \widetilde F_k[u^L](s,\tau) &\doteq
    \widetilde F_k\big(\tau;\,
    u(\,\cdot\ ;\,u^L,s), \,v_k(\,\cdot\ ;\,u^L,s), \,\sigma(\,\cdot\
    ;\,u^L,s)\big)
  \end{array}
  \quad\forall~\tau\in [0,s]\,,
  \label{sigmaFdef}
\end{equation}
for the $\sigma$-component of the solution to (\ref{eq2:Ti}), and for
the reduced flux evaluated in connection with such a solution.

For negative values $s<0$, $|s|<<1,$ one replaces in (\ref{eq2:Ti})
the lower convex envelope of~$\widetilde F_k$ on the interval $[0,s]$
with its upper concave envelope on $[s,0]$
(defined in analogous way as \eqref{convenv}), and then constructs the
curve $T_k[u^L]$ and the map $\sigma_k[u^L]$ exactly in the same way
as above looking at the solution of the integral system (\ref{eq2:Ti})
on the interval $[s,0]$.  The elementary curve $T_k[u^L]$ and the
wave-speed map $\sigma_k[u^L]$ constructed in this way enjoy the
properties stated in  in the following
theorem, where we let $\C_I([a,b])$ denote the set of continuous and
increasing scalar functions defined on an interval $[a,b]$, and we set
$\C_I([a,b])\doteq\C_D([b,a])$ in the case $a>b$, letting
$\C_D([b,a])$ denote the set of continuous and decreasing scalar
functions defined on $[b,a]$.
\begin{theorem}[\cite{BB,srp}]
\label{thm:riemsolv}
Let $A$ be a smooth, matrix valued map defined from a
do\-ma\-in $\Omega\subset\real^N$ into~$\M^{N\times
N}(\real)$, and assume that the matrices $A(u)$ are strictly
hyperbolic.  Then, for every $u\in\Omega$, there exist $N$ Lipschitz
continuous curves $s\to T_k[u](s)\in\Omega$ satisfying
${\displaystyle{\lim_{s\to 0}}}\, \frac{d}{ds}T_k[u](s)=r_k(u)$,
together with $N$ continuous functions $s\to
\sigma_k[u](s,\,\cdot)\in\C_I([0,s])$ $(k=1,\dots,N)$, defined on a
neighborhood of zero, so that the following holds.  Whenever $u^L\in
\Omega,$ $u^R=T_k[u^L](s)$, for some $s$, letting
$\mathcal{I}\doteq\{\tau\in[0,s\,]\,:\,
\sigma_k[u^L](s,\,\tau)\neq\sigma_k[u^L](s,\,\tau')\ \text{for
all}~\tau'\neq \tau\}$, the piecewise continuous function
\begin{equation}
  \label{eq1:solRPi}
  u (t,x) \doteq
  \begin{cases}
    u^L\quad &\text{if\quad $x/t <\sigma_k[u^L](s,\,0)\,,$}\\
    \noalign{\smallskip}
    T_k[u^L] (\tau)\quad &\text{if\quad $x/t = \sigma_k[u^L](s,\,\tau)$
    \ \ \ \ for \ some  \ \ $\tau\in\mathcal{I}$\,,}\\
    \noalign{\smallskip}
    u^R\quad &\text{if\quad $x/t >\sigma_k[u^L] (s,\,s)\,,$}
  \end{cases}
\end{equation}
provides the unique vanishing viscosity solution $($determined by the
parabolic approximation
$(1.16))$ of the Riemann
problem~$(\ref{sys}), (\ref{riemin})$.
\end{theorem}

\begin{remark}
\label{rem2:consform}
If the system (\ref{sys}) is in conservation form, i.e. in the case
where $A(u)=DF(u)$ for some smooth flux function $F$, and if the
characteristic fields satisfy the assumption (H), the general solution
of the Riemann problem provided by (\ref{eq1:solRPi}) is a composed
wave of the $k$-th family made of a finite number of
contact-discontinuities 
(which satisfy the Liu
admissibility condition of Definition~\ref{def:Laxstab}) adjacent to
rarefaction waves. Namely, the regions where the $v_k$-component of
the solution to~(\ref{eq2:Ti}) vanishes correspond to rarefaction
waves if the $\sigma$-component is strictly increasing and to contact
discontinuities if the $\sigma$-component is constant, while the
regions where the $v_k$-component of the solution to~(\ref{eq2:Ti}) is
different from zero correspond to contact discontinuities 
or to compressive shocks.  In particular,
whenever the solution of a Riemann problem with initial data $u^L,\,
u^R=T_k[u^L](s)$ contains a Liu admissible shock joining, say, two
states $T_k[u^L](s'), T_k[u^L](s'')$, $s',s''\in[0,s]$, one has
$\sigma_k[u^L](s,\,s')=\sigma_k[u^L](s,\,\tau)$ for
all~$\tau\in[s',s'']$, and $\sigma_k[u^L](s,\,s')$ provides the shock
speed of the discontinuity $\big(T_k[u^L](s'),\, T_k[u^L](s'')\big)$.
Clearly, in a non conservative setting, ``admissibility'' for a jump
means precisely that the jump corresponds to a traveling profile for
the  parabolic approximation
with identity viscosity matrix~$(1.16)_1$.
\end{remark}

Once we have constructed the elementary curves $T_k$ for each $k$-th
characteristic fa\-mi\-ly, the {\it vanishing viscosity solution} of a
general Riemann problem for (\ref{syscon}) is then obtained by a
standard procedure observing that the composite mapping
\begin{equation}
  (s_1,\dots,s_N)\mapsto T_N(s_N)\circ \cdots \circ T_1(s_1) [u^L]
\doteq u^R\,,
  \label{eq2:RP1}
\end{equation}
is one-to-one from a neighborhood of the origin in $\real^N$ onto a neighborhood of
$u^L$. This is a consequence of the fact that the curves $T_k[u]$ are
tangent to $r_k(u)$ at zero (cfr.~Theorem~\ref{thm:riemsolv}), and
then follows by applying a version of the implicit function theorem
valid for Lipschitz continuous maps.  Therefore, we can uniquely
determine intermediate states $u^L\doteq\omega_0,$ $\omega_1,$
$\dots,$ $\omega_N\doteq u^R$, and wave sizes $s_1, \dots, s_N,$ such
that there holds
\begin{equation}
  \label{eq2:RP2}
  \omega_k =T_k[\omega_{k-1}](s_k)\quad\qquad k=1,\dots,N\,,
\end{equation}
provided that the left and right states $u^L, u^R$ are sufficiently
close to each other.  Each Riemann problem with initial data
\begin{equation}
  \label{elemriem}
  \overline u_k(x) = 
  \begin{cases}
    \omega_{k-1} &\text{if \ \ $x<0$,}\\
    \omega_k &\text{if \ \ $x>0$,}
  \end{cases}
\end{equation}
admits a vanishing viscosity solution of {\it total size} $s_k$,
containing a sequence of rarefactions and Liu admissible
discontinuities of the $k$-th family.  Then, because of the uniform
strict hyperbolicity assumption~(\ref{hyp}), the general solution of
the Riemann Problem with initial data $\big(u^L,\,u^R\big)$ is
obtained by piecing together the vanishing viscosity solutions of the
elementary Riemann problems (\ref{syscon}) (\ref{elemriem}).  Throughout
the paper, with a slight abuse of notation, we shall often call $s$ a
wave of (total) size~$s$, and, if $u^R=T_k[u^L](s)$, we will say that
$(u^L,\,u^R)$ is a wave of size $s$ of the $k$-th characteristic
family.
\medskip

A fundamental ingredient in order to get a convergence rate for the
Glimm scheme is the wave tracing procedure, which was first introduced
by T.P. Liu in his celebrated paper~\cite{tplGlimm} for systems with
genuinely nonlinear or linearly degenerate fields, and lately extended
to systems fulfilling assumption (H)~\cite{tplams,glimmly}. In this spirit, 
we introduce the following notion of partition of a $k$-wave
$(u^L,u^R)$, defined in terms of the elementary
curves $T_k$ at (\ref{Tdef}).
\begin{definition}
\label{def:par}
Given a pair of states $u^L, u^R$, with $u^R = T_k[u^L] (s)$ for some
$s>0$, we say that a set 
$\big\{ y^1, \dots, y^\ell \big\}$ is a
partition of the $k$-th wave $(u^L,u^R)$ 
at time $i\eps$, if the followings holds.
\begin{enumerate}
  \item
    There exist scalars $s^h>0$, $h=1,\dots,l$, such that,
    setting $\tau^h\doteq\sum_{p=1}^h s^p$, $w^h\doteq T_k[u^L](\tau^h)$,
    there holds
    $$
     y^h=w^h-w^{h-1}\qquad\ \ \forall~h\,.
    $$
    The quantity $s^h$ is called the {\it size} of the elementary wave $y^h$.
  \item
    Letting $\sigma\doteq\sigma_k[u^L](s,\cdot)$ be the map in \eqref{sigmaFdef},
    there holds 
    $$
    \sigma(s_h) - \sigma(s_{h-1}) \leq \eps\qquad\ \forall~h\,.
    $$
    Moreover, we require that $\theta_{i+1}\notin\,]\sigma(\tau^{h-1}),\, \sigma(\tau^h)[$\,,
    for all $h$ (so to avoid further partitions of $y^h$ at $t=(i+1)\eps$).
\end{enumerate}
The definition is entirely
similar in the case $u^R = T_k[u^L] (s)$, with $s<0$.
In connection with a partition $\big\{ y^1, \dots, y^\ell \big\}$
of $(u^L,u^R)$, we define the corresponding speed of the elementary wave $y^h$ as
\begin{equation}
\label{eq:parsp}
\lambda_k^h \doteq \frac{1}{s^h}
\int_{\tau^{h-1}}^{\tau^h} \sigma(\tau)~d\tau\qquad\ \ \forall~h\,.
\end{equation}

\end{definition}

%
%
%	ONE INFLECTION POINT
%
%
\section[The case of a single linearly degenerate manifold]{The case of a single linearly degenerate manifold}
\label{sec:oip}
\setcounter{equation}{0}

In this section we will establish the basic estimates on the
change in size and speeds of the elementary waves 
of an approximate
solution provided by the Glimm scheme, under the following simplified
assumption for the hyperbolic system (\ref{syscon}) (or for the quasilinear
sytem \eqref{sys}).
\begin{description}
\item[{\rm (H1)}] 
  For each $k\in\{1,\dots,N\}$-th characteristic family
   the linearly de\-ge\-ne\-rate manifold
  $\mathcal{M}_k$ at (\ref{LDM})
  is either empty (GNL), or it is the whole space (LD), or it consists of \emph{a single} smooth,
  $N\!-\!1$-dimensional, connected, manifold and there holds \eqref{eq:transv-ld} (NGNL).
\end{description}

\noindent
The general solution of a Riemann problem for a sysytem satisfying
the assumption (H1) consists of rarefaction waves,
compressive shock and composed waves made of a single one-side contact discontinuity 
adjacent to a rarefaction wave.
For such systems,
we may consider  the same type of quadratic {\it interaction
potential} introduced in \cite{amwftnpn}  for approximate solutions 
constructed by a front tracking algorithm, which in the case of
solutions $u^\epsilon$ generated by a Glimm scheme can be defined by setting
\begin{equation}
\begin{aligned}
  \label{eq:ip1}
  Q_1 (t) \doteq 2 \sum_{\substack{k_\alpha=k_\beta\\ s_\alpha
      s_\beta>0}} \big\vert s_\alpha
  s_\beta \big\vert &+ 2 \, \sum_\alpha \big\vert
  s_\alpha^r s_\alpha^s \big\vert + \sum_\alpha \big\vert
  s_\alpha^r \big\vert^2 +
 \\
  &+ \,c_0 \Bigg[ \sum_{\substack{k_\alpha=k_\beta\\ s_\alpha
  s_\beta<0}} 
  + \sum_{\substack{k_\alpha<k_\beta\\ x_\alpha(t)>x_\beta(t)}}
   \Bigg]\big\vert s_\alpha s_\beta \big\vert\,,
\end{aligned}
\end{equation}
where $c_0>2$ is a suitable large constant to be defined later, 
$s_\alpha$ denotes the size of a wave of the $k_\alpha$-th family of $u^\epsilon(t)$ located at $x_\alpha(t)$,
while $s_\alpha^r$, $s_\alpha^s$ are, respectively, the (possibly zero) rarefaction and shock
components of a  wave~$s_\alpha$.
The presence of the factor 2 in the first two summands guarantees
the invariance of $Q_1$ when two portions of rarefaction fans  
of the same family,
emanating from two consecutive mesh-points,
are joined together for the effect of sampling, since otherwise the
quantity $Q_1$ would increase for the presence of the square 
of the rarefaction components.
As customary, we shall define the {\it total strength
of waves} in $u^\epsilon(t)$ as
\begin{equation}
 \label{eq:v}
V (t) \doteq \sum_\alpha |s_\alpha|\,.
\end{equation}
\smallskip

To fix the ideas, assume that the second derivative of
$\lambda_k$ in \eqref{eq:transv-ld} is negative, i.e. that
\begin{equation}
  \label{eq:transv-ld1}
   \nabla (\nabla\lambda_k \cdot r_k)(u)\cdot r_k(u)< 0\qquad\forall u\in \mathcal{M}_k\,.
\end{equation}
In order to control the nonlinear coupling of waves of the
same family and with the same sign of two Riemann solutions
for sysytems satisfying
the assumption (H1),
 as in \cite{amwftnpn} we introduce the following
definition of  quantity of interaction.
\begin{definition}
\label{def:qi-sld}
Consider two nearby waves
of sizes $s', s''$ with the same sign and belonging to  the the same $k$-th characteristic family, 
with $s'$ located at the left of $s''$.
Let $u', u''$ be the left state of $s', s''$, respectively,
and assume that there exist waves $s'_i$, $k<i\leq N$, of the $i$-th family,
$s''_j$, $1\leq j<k$, of the $j$-th family,
so that $u''=\big(\bigcirc_{j=1}^{k-1} T_j(s''_j)\big)\circ\big(\bigcirc_{i=k}^N T_i(s'_i)\big)[u']$.
Then, we define the {\it quantity of interaction} between $s'$ and $s''$ as
\begin{equation}
  \label{eq:qi1}
  {I}_1 (s', s'')\doteq
 \big( \big\vert (s' + s'')^r - s'^r \big\vert +
    \vert s'^s \vert \big) \vert s''^s \vert\,,
\end{equation}
where $s'+s''$ must be interpreted as the size of a $k$-wave
having left state $u^\sharp\doteq \bigcirc_{j=1}^{k-1} T_j(s''_j)[u']$,
while $s^r$, $s^s$ denote, respectively, the (possibly zero) rarefaction and shock
components of a wave $s$.
\end{definition}
\begin{remark}
\label{rem2b:qi-sld}
In the case where $s', s''$ are both rarefactions the quantity of interaction $I_1$ 
in \eqref{eq:qi1} vanishes, while $I_1(s', s'')=|s' s''|$ whenever $s' s''$
are both shock waves.
\end{remark}
\medskip

By standard arguments (e.g, see \cite[Section~9.6, Section~13.4]{Daf})
one can obtain as in \cite{amwftnpn}
the basic estimates on the change in values of
the total strength of waves $V(t)$  and of the
interaction potential $Q_1(t)$, across the grid-times $i \eps$,
for an approximate solution $u^\epsilon$ constructed by the Glimm scheme.
Namely, defining for every pair of waves of the same family $s', s''$
the {\it amount of cancellation} $\C(s', s'')$ as
\begin{equation}
\label{eq:can}
\C(s',s'')\doteq
\begin{cases}
  \min \big\{ \vert s'\vert, \vert s''\vert \big\}\quad
  &\text{if\quad $s' s''<0$\,,}\\
  \noalign{\smallskip}
  0\quad & \text{otherwise,}
\end{cases}
\end{equation}
the
following generalization of \cite[Lemma~2.1, Lemma~5.1]{amwftnpn} hold.
\begin{lemma}
\label{lem:ie1}
Under the assumption (H1),
let 
$s'_1,\ldots,s'_N$ and $s''_1,\ldots,s''_N$ be, respectively, the sizes of the
waves in the solution of two adjacent Riemann problems $(u^L, u^M)$
and $(u^M, u^R)$,  $s'_i$
and $s''_i$ belonging to the $i$-th characteristic family. Call
$s_1,\ldots,s_N$ the sizes of the
waves in the solution of the Riemann problem $(u^L, u^R)$, 
$s_i$ belonging to the $i$-th characteristic family. Then, there holds
\begin{equation}
\label{eq:ie1}
\sum_{k=1}^N \big\vert s_k-s'_k-s''_k \big\vert = \mathcal{O}(1) \cdot \Bigg[
  \sum_{\substack{1\leq i,j \leq  N\\i> j}} \vert s'_i s''_j \vert +
 \!\!\!\sum_{\substack{i=1,\ldots,N\\ s'_i s''_i<0}}
  \vert s'_i s''_i \vert +
 \!\!\!\sum_{\substack{i=1,\ldots,N\\ s'_i s''_i>0}}
   {I}_1 (s'_i, s''_i) \Bigg]\,.
\end{equation}
Moreover, for any $k\in \{ 1,\ldots,N \}$-th NGNL characteristic family,  the following estimates on the rarefaction components of
the outgoing waves hold.
\begin{align}
\null
&\big\vert s_k^r - (s'_k+s''_k)^r \big\vert = \mathcal{O}(1) \cdot\! \Bigg[
  \sum_{\substack{1\leq i,j \leq  N\\i> j}} \vert s'_i s''_j \vert +
   {I}_1 (s'_k, s''_k) \Bigg]\nonumber\\
\label{eq:rar11}
&\hspace{7.5truecm} \text{if}\quad s'_k s''_k>0\,,\\
\noalign{\medskip}
&\big\vert s_k^r - (s'^r_k+s''^r_k) \big\vert = 
\mathcal{O}(1) \cdot \Bigg[
  \sum_{\substack{1\leq i,j \leq N\\i> j}} \vert s'_i s''_j \vert +
 \min\{|s'_k|,|s''_k|\}  \Bigg]\nonumber\\
\label{eq:rar12}
&\hspace{7.5truecm} \text{if}\quad s'_k s''_k<0\,,
\end{align}
where, in
(\ref{eq:rar11}) $s'_k+s''_k$ represents the size of a $k$-wave having
left state $u^\sharp\doteq \big(\bigcirc_{j=1}^{k-1} T_j(s''_j)\big)\circ\big(\bigcirc_{i=1}^{k-1} T_i(s'_i)\big)[u^L]$.
\end{lemma}
\begin{lemma}
\label{lem:ie12}
In the same setting of Lemma~\ref{lem:ie1},
provided that the total strength of waves is sufficiently small,
there exists some constant $c_0>0$ (in \eqref{eq:ip1})
so that there holds
\begin{align}
  \label{eq:tv1est}
  \Delta V  &\leq-\sum_{1\leq i \leq N} \C(s'_i,s''_i) + \mathcal{O}(1) \cdot \Bigg[
  \sum_{\substack{1\leq i,j \leq  N\\  \noalign{\smallskip} i> j}} \vert s'_i s''_j \vert +
 \!\!\!\sum_{\substack{1\leq i \leq N\\  \noalign{\smallskip}s'_i s''_i>0}}
   {I}_1 (s'_i, s''_i) \Bigg]\,,\\
  \noalign{\smallskip}
  \label{eq:ip1est}
  \Delta Q_1  &\leq -\frac{1}{2} \Bigg[
  \sum_{\substack{1\leq i,j \leq  N\\  \noalign{\smallskip} i> j}} \vert s'_i s''_j \vert +
 \!\!\!\sum_{\substack{1\leq i \leq N\\ \noalign{\smallskip} s'_i s''_i<0}}
  \vert s'_i s''_i \vert +
 \!\!\!\sum_{\substack{1\leq i \leq N\\  \noalign{\smallskip} s'_i s''_i>0}}
   {I}_1 (s'_i, s''_i) \Bigg]\,.
\end{align}
Here, as customary, we use the notations 
$\Delta V\doteq V^+ - V^-$, $\Delta Q_1 \doteq Q_1^+ - Q_1^-$,
where $V^-, Q_1^-$ and $V^+, Q_1^+$ denote, 
respectively, the values of
$V, Q_1$ related to the incoming waves $s'_1,\ldots,s'_N$, $s''_1,\ldots,s''_N$,
and to the outgoing waves $s_1,\ldots,s_N$.
\end{lemma}
Relying on Lemma~\ref{lem:ie12}, one deduces that there exists some 
constant $C_1>0$, independent of $\eps$, so that
if $V(t), Q_1(t)$ denote 
the total strength of waves and  the
interaction potential of an approximate solution $u^\epsilon(t)$ 
constructed by the Glimm scheme, 
the functional
\begin{equation}
  \label{eq:ups1}
  t\mapsto \ups_1(t) \doteq V(t) + C_1 Q_1 (t)
\end{equation}
is non increasing at any time, provided that the total initial 
strength $V(0)$
is sufficiently small.
Moreover, for any given $0\leq m<n$, the total amount of wave interaction and cancellation 
taking place in
the time interval $[m\eps, n\eps]$ is bounded by $\mathcal{O}(1)\cdot
(\ups_1 (m\eps) - \ups_1 (n\eps))$.
Denote $\Delta \ups^{m,n}_1 \doteq \ups_1 (n\eps) - \ups_1 (m\eps)$
 the variation of $\ups_1$ on $[m\eps, n\eps]$.
\medskip

A basic ingredient of the strategy followed in \cite{bm} to
establish a convergence rate of the Glimm scheme 
is the wave tracing algorithm introduced in~\cite{tplams} 
for GNL or LD systems, 
and then extended in \cite{glimmly} to NGNL systems,
which consists in partioning
the outgoing waves issuing from every mesh point $(i\eps, j\eps)$ 
in two type of waves: \emph{primary waves} (i.e. waves that can be
traced back from the time $t=i\eps$ to a previous time $t=m\eps<i\eps$), and
\emph{secondary waves} (i.e. waves that are generated by
interactions occurring in the time interval  $]m\eps,i\eps]$, or 
that are canceled before a later time $t=n\eps>i\eps$).
The total strength of secondary waves produced in a given time interval 
$[m\eps,n\eps]$ is bounded by the total amount of interaction 
and cancellation occurring within  $[m\eps,n\eps]$.

The key step of this procedure is to show that 
the variation of a Glimm functional
provides a bound for
the change in strength and for the product of strength 
times the variation in speeds of the primary waves.
The main novelty of the analysis performed here 
consists in implementing a wave tracing algorithm
for a NGNL system satisfying the assumption~(H1)
in which such bounds are obtained
relying on a Glimm functional with a quadratic
potential interaction, differently from the Glimm functional
with a cubic potential interaction  used in \cite{glimmly}.
Namely, 
recalling the
Definition~\ref{def:par} of a wave partition, we have the following result.
\begin{proposition}
\label{pro:wt1}
Under the assumption (H1),
given a Glimm approximate solution and any  fixed $0\leq m<n$, 
there exists a partition of elementary wave
sizes and speeds
$\big\{ y^h_k (i,j),\, \lambda_k^h(i,j) \big\}$, $k=1,\dots, N$,
$i=m,m+1,\dots,n$, $j\in\integer$, so that the following hold.
\begin{enumerate}
\item
For every $i, j, k$, 
$\big\{ y^h_k (i,j) \big\}_{0<h\leq \ell_k(i,j)}$ is a partition of
the wave of the $k$-th family issuing from $(i\eps,j\eps)$, and $\big\{
\lambda_k^h(i,j) \big\}_{0<h\leq \ell_k(i,j)}$ are the corresponding speeds, 
according with 
Definition~\ref{def:par}.
\item
For every $i, j, k$, 
$\big\{ y^h_k (i,j),\, \lambda_k^h(i,j) \big\}_{0<h\leq \ell_k(i,j)}$ is a disjoint union of
the two sets
$$
\big\{ \wy^h_k (i,j),\, \wl_k^h(i,j) \big\}\,, \quad   \quad
\big\{ \wwyhk (i,j),\, \wwl_k^h (i,j) \big\}\,,
$$
with the following properties:
\begin{enumerate}
\item
\begin{equation}
\label{eq:ttest1}
\sum_{j,k,h} \big\Vert \wwyhk (i,j) \big\Vert = \mathcal{O}(1) \cdot
\big|\Delta \ups_1^{m,n}\big|
\qquad\quad\forall~m\leq i \leq n\,;
\end{equation}
\item
for every fixed $i, k, h$, there is a one-to-one correspondence between 
$\big\{ \wy^h_k (m,j),\wl_k^h(m,j)~;~ j\in\integer \big\}$ and \,$\big\{ \wy^h_k (i,j), \wl_k^h(i,j)~;~ j\in\integer \big\}$:
\begin{equation}
\label{eq:tcor1}
\big\{ \wy^h_k (m,j), \wl_k^h(m,j) \big\} \leftrightarrow
\big\{ \wy^h_k (i,\ell_{(i,j,k,h)}), \wl_k^h(i,\ell_{(i,j,k,h)}) \big\}
\end{equation}
such that the sizes $\ws^h_k$ and the speeds $\wl_k^h$ of
the corresponding waves satisfy
\begin{gather}
\label{eq:stest1}
\sum_{j,k,h} \left( \max_{m\leq i\leq n} \big\vert \ws^h_k (m,j) -
\ws^h_k (i,\ell_{(i,j,k,h)}) \big\vert \right) = \mathcal{O}(1)
\cdot \big|\Delta \ups_1^{m,n}\big|\,,\\
\noalign{\smallskip}
\label{eq:spest1}
\!\!\!\!\!\!\!\!\!
\sum_{j,k,h} \left( \big\vert \ws^h_k (m,j) \big\vert \cdot \!
\max_{m\leq i\leq n} \big\vert \wl^h_k (m,j) -
\wl^h_k (i,\ell_{(i,j,k,h)}) \big\vert \right) = \mathcal{O}(1)
\cdot \big|\Delta \ups_1^{m,n}\big|\,.
\end{gather}
\end{enumerate}
\end{enumerate}
\end{proposition}
\textbf{Proof.}
The desired partition for
an approximate solution $u^\eps$ will be
constructed proceeding 
by induction on the time steps $i\eps$, $m\leq i \leq n$.
Assuming that a partition of elementary waves
fulfilling properties 1-2 is given for all times $m\eps\leq t< i\eps$,
we wish to produce a partition of the outgoing waves generated by the
interactions occurring at $t=i\eps$, so to preserve the pro\-per\-ties~1-2.
Observe first that the existence of such a partition is
already guaranteed by the analysis in \cite{glimmly} 
if all interactions take place between waves of different family or of the
same family with opposite sign, since for systems satisfying the
assumption (H1)
the change in strength and the product of strength 
times the variation in speeds of the primary waves 
is controlled by the variation of a Glimm functional with
quadratic interaction potential as the
part in brackets of \eqref{eq:ip1}.

Therefore, it will be sufficient to consider an
interaction between two waves 
issuing from two consecutive mesh points $((i-1)\eps, (j-1)\eps)$ and $((i-1)\eps,
j\eps)$,
say $s_k'$, $s_k''$, belonging to a $k$-th NGNL characteristic family,
and having the same sign.
For the sake of simplicity, assume that $s_k',s_k''>0$.
Let $s_p$ ($p=1,\dots,N$) be the outgoing wave 
of the $p$-th family 
issuing from $(i\eps, j\eps)$,
and let
\begin{equation}
\label{partinc}
\big\{ y'^h_k,\, \lambda'^h_k \big\}_{0<h\leq \ell'}\,,
\qquad\quad 
\big\{ y''^h_k,\, \lambda''^h_k \big\}_{0<h\leq \ell''}\,, 
\end{equation}
be the partitions of $s'_k$ and $s''_k$ enjoing the properties 1-2
on $[m\eps, (i-1)\eps]$, with sizes
\begin{equation}
  \big\{s'^h_k \big\}_{0<h\leq \ell'}\,,
 \quad\quad 
\big\{s''^h_k\big\}_{0<h\leq \ell''}\,.
\end{equation}
For every $p\neq k$-th wave $s_p$, we may choose a partition
$\{ y^h_p \}_{0<h\leq \ell_p}$ as in Definition~\ref{def:par}, 
with corresponding speeds $\{ \lambda^h_p \}_{0<h\leq \ell_p}$.
Then, if we label all the subwaves $y^h_p$ as secondary waves 
${\widetilde{\widetilde y}}^h_p$,
the bound \eqref{eq:ttest1} (for $i,j,p$) is certainly satisfied thanks to the
interaction estimates \eqref{eq:ie1}.
Instead, for the $k$-th wave $s_k$, possibly considering a refinement
of the partition of $s''_k$ (or of $s'_k$)
we may assume that either 
$s'_k+s''^1_k\leq s_k$, or $s'^1_k\leq s_k$ 
(in the case  $s'_k\geq s_k$),
and let $\overline \ell'\doteq \max\{h\leq \ell'\,:\, \sum_{q=1}^h s'^q_k
\leq s_k\}$, $\overline \ell''\doteq \max\{h\leq \ell''\,:\, s'_k+\sum_{q=1}^h s''^q_k
\leq s_k\}$.
Then, we define a partition of $s_k$ by means of its sizes, 
setting
\begin{equation}
  \label{eq:par12}
  s^h_k \doteq
  \begin{cases}
    s'^h_k\quad &\text{if\quad
      $h=1,\ldots,\overline\ell'$\,,}\\
    \noalign{\medskip}
    s''^{h-\ell'}_k\quad &\text{if\quad $\overline\ell'=\ell'$\quad and\quad
      $h=\ell'+1,\ldots,\ell'+\overline \ell''$}\\
  \end{cases}
\end{equation}
(possibly refining the partitions \eqref{partinc} so to satisfy 
property 2 of Definition~\ref{def:par}),
and choosing a partition of $s_k-(s'_k+s''_k)$
as in Definition~\ref{def:par} in the case $s_k > s'_k+s''_k$.
The subwaves $s^h_k$ in \eqref{eq:par12} inherit the same classification
in primary and secondary waves of the corresponding subwaves
$s'^h_k$ or $s''^{h-\ell'}_k$, while all the possible subwaves of 
$s_k-(s'_k+s''_k)$ are labelled as secondary waves.
Clearly, the bound \eqref{eq:ttest1} is again satisfied because of the
interaction estimates \eqref{eq:ie1}, while the
ono-to-one correspondence at (\ref{eq:tcor1}) and the bound 
\eqref{eq:stest1} are verified by construction and by the inductive assumption.
Hence, in order to conclude the proof, it remains to establish
only the estimate \eqref{eq:spest1}.

By the assumption (H1), and because the incoming waves $s'_k, s''_k$
have the same sign, at most one of them can possibly be a 
composed wave, say~$s'_k$, while $s''_k$ will be a shock.
Denote as $(s'_k)^r, \, (s'_k)^s$ the rarefaction and shockcomponent
of $s'_k$, respectively.
For sake of simplicity, assume that $s_k>s'_k$, i.e. that $\overline \ell'=\ell'$.
The outgoing wave $s_k$ is either a shock or a composed wave.
In the first case its Rankine-Hugoniot speed $\lambda_k$ coincides with
the speeds $\lambda_k^h$ of all subwaves $s^h_k$ defined according with 
Definition~\ref{def:par}, since for
a shock wave the integrand function in \eqref{eq:parsp}
results to be a constant (cfr. Remark~\ref{rem2:consform}).
 Hence, letting $(\lambda'_{k})^s,\, \lambda''_k$ denote
the speeds of the shock component of $s'_k$ and of $s''_k$,
respectively, by a direct computation one finds
\begin{equation}
\begin{aligned}
  &\big\vert \lambda^h_k - \lambda'^h_k \big\vert \leq 
    \big\vert \lambda_k - (\lambda'_{k})^s \big\vert= \mathcal{O}(1)\cdot
   s''_k&\quad &\forall~h=1,\ldots,\ell'\,,&\\
  &\big\vert \lambda^h_k - \lambda''^h_k \big\vert =
    \big\vert \lambda_k - \lambda''_{k} \big\vert =
  \mathcal{O}(1)\cdot  s'_k&\quad &\forall~
  h=\ell'+1,\ldots,\ell'+\overline \ell''\,.&
\end{aligned}
\label{eq:est1s4}
\end{equation}
In turn, \eqref{eq:est1s4} implies
\begin{equation}
  \sum_{h=1}^{\ell'} s'^h_k\big\vert \lambda^h_k - \lambda'^h_k \big\vert 
  +\sum_{h=\ell'+1}^{\ell'+\overline\ell''} s''^h_k
\big\vert \lambda^h_k - \lambda''^h_k \big\vert
=  \mathcal{O}(1)\cdot  s'_k s''_k\,,
\label{eq:ests1s5}
\end{equation}
which, relying on the inductive assumption, yields \eqref{eq:spest1} since in this case,
by the estimate \eqref{eq:rar11}, and because the rarefaction component of $s_k$
is zero, there holds
$s'_k s''_k  = \mathcal{O}(1) \cdot{I}_1 (s'_k,s''_k)$.

Next, assume that the outgoing wave $s_k$ is made of a rarefaction component
$(s_k)^r<(s'_k)^r$ and of a shock component $(s_k)^s$. 
Then, possibly considering a refinement
of the partition of $s_k$ (and hence of the partition of~$s'_k$), 
there will be some index $\ell^r<\ell'$ so that $\sum_{h=1}^{\ell^r} s_k^h =(s_k)^r$
and $\sum_{h=\ell^r+1}^{\ell'+\overline\ell''} s_k^h =(s_k)^s$.
Notice that $(s_k)^s$ can be seen as a
shock wave generated by an interaction between a composed wave
with rarefaction component of size $(s'_k)^r-(s_k)^r$ and 
shock component $(s'_k)^s$, and of the shock wave $s''_k$,
for which we can apply the above estimates on the variation of wave speeds.
Hence, the wave speeds $\lambda_k^h$ of $s_k^h$ 
defined as in \eqref{eq:parsp} satisfy
\begin{equation}
\begin{gathered}
\big\vert \lambda_k^h - \lambda'^h_k \big\vert =
\begin{cases}
0\quad &\forall~h=1,\ldots,\ell^r\,,\\
\noalign{\smallskip}
\mathcal{O}(1) \cdot s''_k\quad
&\forall~h=h^r,\ldots,\ell'\,,
\end{cases}
\\
\noalign{\smallskip}
\big\vert \lambda_k^h - \lambda''^h_k \big\vert =
\mathcal{O}(1)\cdot  \big[(s'_k)^r-(s_k)^r+(s'_k)^s\big]\qquad\ \ \forall~
  h=\ell'+1,\ldots,\ell'+\overline \ell''\,.
\end{gathered}
\label{eq:est1s6}
\end{equation}
This implies
\begin{equation}
  \sum_{h=1}^{\ell'} s'^h_k\big\vert \lambda^h_k - \lambda'^h_k \big\vert 
  +\sum_{h=\ell'+1}^{\ell'+\overline\ell''} s''^h_k
\big\vert \lambda^h_k - \lambda''^h_k \big\vert
=  \mathcal{O}(1)\cdot  \big[((s'_k)^r-(s_k)^r)+(s'_k)^s\big]\, s''_k\,,
\label{eq:ests1s7}
\end{equation}
which in turn, relying on the inductive assumption, yields again \eqref{eq:spest1} since in this case,
by the estimate \eqref{eq:rar11} and because $s''_k=(s''_k)^s$, one has
$\big[(s'_k)^r-(s_k)^r+(s'_k)^s\big]\, s''_k  = \mathcal{O}(1) \cdot{I}_1 (s'_k,s''_k)$.
 This completes
the proof of the proposition.
\qed

Proposition~\ref{pro:wt1} provides for NGNL systems
satisfying the assumption (H1) the same type of result 
that was established in \cite[Proposition~2]{bm} for systems
with GNL or LD characteristic families. In order to
obtain the desired convergence rate \eqref{eq:Glimm}
one can now simply repeat the proofs of \cite[Propositions~3-4]{bm}
and of the final estimates in \cite[\S~6]{bm}, which all rely 
only on the conclusion of \cite[Proposition~2]{bm} and thus remain
valid within our more general framework of NGNL systems.
We will give a brief description of them in
Section~\ref{sec:con}.

\begin{remark}
\label{rem3:c11ngnl}
The conclusion of Theorem~\ref{thm:glimm}, established so far for
smooth systems satisfying the assumption (H1), remains valid 
if we assume that the flux function $F$ is $\C^{2,1}$
and that, for each $k$-th characteristic family not fulfillying (H1),
the linearly degenerate manifold $\MC_k$ in \eqref{LDM} is a $\C^{1,1}$
$N$-dimensional, connected manifold, $F$ is $\C^3$ on $\Omega\setminus\MC_k$,
the vector field $r_k$ is transversal to the boundary of $\MC_k$,
and letting $\partial^+ \MC_k,\, \partial^- \MC_k$ denote the connected components of
the boundary of $\MC_k$
where $r_k$ points towards $\Omega\setminus\MC_k$ and~$\MC_k$, respectively,
there holds
\begin{equation}
  \label{eq:transv-ld2}
\begin{aligned}
   &\nabla^+ (\nabla\lambda_k \cdot r_k)(u)\cdot r_k(u)< 0 
\qquad \quad \forall u\in \partial^+\mathcal{M}_{k}\,,
\\
   &\nabla^- (\nabla\lambda_k \cdot r_k)(u)\cdot r_k(u)< 0 
\qquad \quad \forall u\in \partial^-\mathcal{M}_{k}
\end{aligned}
\end{equation}
($\nabla^{\pm} (\nabla\lambda_k \cdot r_k)(u)\cdot r_k(u)\doteq
{\displaystyle{\lim_{h\to 0\pm}}}\frac{\nabla\lambda_k \cdot 
r_k(u+h r_k(u))-\nabla\lambda_k \cdot r_k(u)}{h}$
denoting the one-side second derivatives of $\lambda_k$).
Indeed, the only difference in the structure of the
elementary waves of a NGNL $k$-th family satisfying such assumptions
instead of (H1) comes from the possible presence of two-sided contact discontinuities.
In fact, under the above assumptions,
the general solution of a Riemann problem of the $k$-th family will be
either a rarefaction wave, or a shock wave (which can be either
a compressive shock or a  contact discontinuity), or a composed wave made of
a rarefaction wave adjacent to one (one-sided
or two-sided) contact discontinuity or several (two-sided)
contact discontinuities.
Then, we may consider the interaction potential $Q_1$ in~\eqref{eq:ip1},
where the shock component of a composed wave $s_\alpha$ containing
several contact discontinuities $s_{\alpha,1}^s, \dots ,  s_{\alpha,l}^s$
is $s_\alpha^s\doteq \sum_{p=1}^l s_{\alpha,p}^s$, and for every such wave 
we add the term $2\,\sum_{\substack{p,q\\ p\neq q}} |s_{\alpha,p}^s\,s_{\alpha,q}^s|$.
One can easily verify that employing this definition of $Q_1$ 
and the same definition of quantity of interaction~$I_1$ in  \eqref{eq:qi1},
the estimates stated in Lemma~\ref{lem:ie1} continue to hold,
provided that
\begin{equation}
  \label{eq:ldml}
\inf\Big\{s>0\,:\,R_k[u](s)\in\partial^+ \MC_k\,,\ u\in\partial^- \MC_k\Big\}>0
\end{equation}
($R_k[u](s)$ denoting the integral curve of $r_k$),
which is certainly true up to a possible slight restriction of the domain $\Omega$.
Relying on Lemma~\ref{lem:ie1}, one then deduces Lemma~\ref{lem:ie12}
and thus can establish the key Proposition~\ref{pro:wt1} with the same arguments
as above.
\end{remark}
%
%
%	INTERACTION POTENTIAL
%
%
\section[A new interaction potential]{A new interaction potential}
\label{sec:ip}
\setcounter{equation}{0}

We turn now our attention to an approximate solution $u^\eps$
constructed by the Glimm scheme for an hyperbolic system
(\ref{syscon}) (or for the quasilinear system \eqref{sys})
that satisfy the 
assumption (H) stated in the Introduction. 
We recall that for such systems
the general solution of a Riemann problem contains
composed waves
made of several contact discontinuities adjacent to
rarefaction waves (instead of just a single contact
discontinuity adjacent to a rarefaction wave as for
the systems treated in \S~\ref{sec:oip}).
We will say that a wave $s$ of this type, belonging to the
$k$-th characteristic family, {\it crosses} all connected components of $\MC_k$
that are transversal to the $k$-th elementary curve $T_k$ issuing from the left state of $s$
and terminating on the right state of $s$.
Notice that, for each $k$-th NGNL family, and for every 
connected component~$\mathcal{M}_{k,h}$ of~$\mathcal{M}_k$,
the first derivative
$\nabla\lambda_k\cdot r_k$ has opposite signs on the  connected
components of $\Omega\setminus\mathcal{M}_{k,h}$ 
adjacent to $\mathcal{M}_{k,h}$, and as a consequence  the 
second derivative $\nabla(\nabla \lambda_k\cdot r_k)\cdot r_k$ has opposite signs on any pair 
of consecutive components $\mathcal{M}_{k,h}, \mathcal{M}_{k,h+1}$.
Thus, by continuity we may assume that there exists some constant
$\delta_0>0$ so that
\begin{equation}
  \label{eq:delta0}
  \nabla(\nabla \lambda_k\cdot r_k)(u)\cdot r_k(u)  \neq 0
   \qquad\ \forall~u \ \ \text{s.t.} \ \ d(u,\, \mathcal{M}_k)\leq 6\delta_0\,,
\end{equation}
where $d(u,\, \mathcal{M}_k)\doteq\inf_{w\in\mathcal{M}_k} |u-w|$
denotes the distance of a state $u$ from~$\mathcal{M}_k$.
\begin{remark}
\label{rem:comld0}
Condition\eqref{eq:delta0} implies that every wave $s$ of a $k$-th NGNL family
with strength $\vert s \vert \leq 3\delta_0$ 
crosses at most one connected component of~$\mathcal{M}_k$.
Moreover,
if an interaction takes place between two waves
of the $k$-th characteristic family  with strength
$\leq \delta_0$, then, by  the interaction estimates in
\cite[Theorem~3.7]{sie}, the outgoing wave of the $k$-th family 
crosses as well at most one connected component of $\mathcal{M}_k$.
\end{remark}

By Remark~\ref{rem:comld0}, as far as the waves of the NGNL families involved in an
interaction have all strength smaller than $\delta_0$, we can 
establish the same kind of estimates of Proposition~\ref{pro:wt1}
employing the quadratic interaction potential in \eqref{eq:ip1} even
for systems satisfying the more general assumption~(H).
On the other hand, observe that if we consider 
an interaction between two shock waves of a $k$-th NGNL family,
say $s',\, s''$, with speeds $\lambda',\, \lambda''$, respectively,
and we assume that $s',\, s''$ have the same sign, then, letting
$\lambda$ denote the shock speed of the outgoing wave 
of the $k$-th family, by the interaction estimates in 
\cite[Theorem~3.7]{sie} there holds
\begin{equation}
  \label{eq:varlex}
[s\Delta\lambda]\doteq
|s'| \big|\lambda-\lambda'\big|+|s''| \big|\lambda-\lambda''\big|
=\OL(1)\cdot\frac{\big|s's''\big|\big|\lambda'-\lambda''\big|}{\big|s'+s''\big|}\,.
\end{equation}
Notice that $\big|s's''\big|\big|\lambda'-\lambda''\big|$ has
precisely the same order of the quantity of which it decreases the interaction potential $\Q$
introduced in \cite{sie} whenever interactions
of this type take place. Therefore, if we assume that at least
one of the incoming waves of the $k$-th family has strength $\geq\delta_0$, we deduce 
from \eqref{eq:varlex} that
$[s\Delta\lambda]=\OL(1)\cdot |\Delta\Q|/\delta_0$.
Hence, for such interactions one may derive the same kind
of estimates on the products of the wave strengths times
the variation of the wave speeds of Proposition~\ref{pro:wt1}
employing the cubic interaction potential  $\Q$
defined in \cite{sie}.

In view of the above observations, we shall introduce now 
a functional $Q$ that is the sum of a quadratic
and of a cubic interaction potential. 
The latter is the interaction potential
for waves of the same family and with the same sign
defined in~\cite{sie}, valid for general strictly hyperbolic systems \eqref{sys},
which takes the form
\begin{equation}
\label{eq:ipngnl}
\Q(t)\doteq \sum_{\substack{k_\alpha=k_\beta\\ s_{\alpha}s_{\beta}>0}} \left\vert \int_0^{s_{\alpha}}
  \int_0^{s_{\beta}} \big\vert \sigma_{\alpha} (\tau) -
  \sigma_{\beta} (\tau') \big\vert~ d\tau d\tau' \right\vert\,.
\end{equation}
The summation here extends to 
all pair of waves $s_\alpha,\, s_\beta$  of the $k_\alpha\in\{1,\dots,N\}$
family with the same sign (including $s_\alpha=s_\beta$), of the approximate solution $u^\eps(t)$,
and $\sigma_\alpha\doteq\sigma_{k_\alpha}[\omega_{\alpha}](s_\alpha,\,\cdot)$ denotes 
the map in \eqref{sigmaFdef}, where $\omega_{\alpha}$ is the left state of $s_\alpha$.
Such a functional controls the nonlinear coupling of waves of the
same family with the same sign.

The quadratic part $Q_q$ of the
functional $Q$ enjoys two basic properties:
\begin{enumerate}
  \item
    it decreases whenever it takes place an interactions between ``small'' waves
    of the same family, i.e. waves
    whose strength is smaller than $\delta_0$, and the amount of
    decreasing satisfies the same type of estimate (\ref{eq:ip1est})
    obtained for
     systems with a single linearly degenerate manifold;
  \item
     the possible
    increase of $Q_q$ caused by interactions involving ``large'' waves of the
   same family, i.e. waves of strength  larger than
   $\delta_0$,
    is controlled by the decrease of $\Q$.
\end{enumerate}
Thus, for general hyperbolic systems \eqref{sys}
satisfying the assumption~(H), we shall consider a potential interaction
of the form
\begin{equation}
  \label{eq:ip}
  Q(t) \doteq Q_q(t) + c\, \Q (t)\,,
\end{equation}
where $c>2$ is a suitable  constant to be specified later.
\smallskip

Towards the defintion of $Q_q$, let us first introduce some further
notations. Given a composed wave $s$ of a $k$-th NGNL family, let $\{ s^h \}_{h =1,\ldots,l}$
be its decomposition in rarefaction and shock components, and write
$h\in \mathcal{R}$ (respectively  $h\in\mathcal{S}$) if $s^h$ is a rarefaction
(respectively a shock) wave. Thus, letting $w^{h-1}, w^h$ denote the left and
right states of each wave $s^h$, one has $w^h=T_k[w^{h-1}](s^h)$.
Next, for every given shock 
$s^h$, $h\in\mathcal{S}$, we define a {\it convex-concave} sub decomposition $\{ s^{h,p} \}_{p =1,\ldots,q_h}$
as follows. Assuming for the sake of simplicity
that $s^h>0$, let $0=\tau^0 <\tau^1 <\cdots < \tau^{q_h}=s^h$
be a partition of $s^h$ determined by the inflection points of
the reduced flux 
$\tau\mapsto\widetilde f_k^h(\tau)\doteq\widetilde F_k[w^{h-1}](s^h,\tau)$ in~\eqref{sigmaFdef},
and set $s^{h,p}\doteq \tau^p-\tau^{p-1}$. We will write $p\in\,\smile$
(respectively $p\in\,\frown$) if $\widetilde f_k^h$ is  convex (respectively  concave) 
on~$[\tau^{p-1},\tau^p]$, and we will call $s^{h,p}$ a {\it convex} (respectively {\it concave})
{\it component} of $s^h$ if $p\in\,\smile$
(respectively $p\in\,\frown$). Then, considering the affine map
\begin{equation}
  \label{eq:vfi}
  \vfi (s)\doteq
  \begin{cases}
    1\quad &\text{if\quad $\vert s\vert\geq 2\delta_0$\,,}\\
    \noalign{\smallskip}
     \big( \vert s\vert-\delta_0 \big)/\delta_0\quad
    &\text{if\quad $\delta_0 \leq \vert s\vert < 2\delta_0$\,,}\\
    \noalign{\smallskip}
    0\quad &\text{if \quad $\vert s\vert \leq \delta_0$\,,}
  \end{cases}
\end{equation}
we define
the {\it intrinsic interaction potential} of $s^h$, $h\in \SC$, as
\begin{equation}
  \label{eq:iip1}
  q(s^h) = \vfi(s^h) \cdot \Bigg[ 2\,\sum_{p\neq q}
  \big\vert s^{h,p} s^{h,q} \big\vert +
  \sum_{p\in\,\smile} \big\vert s^{h,p} \big\vert^2 \Bigg]\,,
\end{equation}
where the first summand runs over all indexes $p,q\in \,\smile\,\cup\,\frown$,
$p\neq q$, and 
$q(s^h)$ is understood to be zero if $s^h$ has zero convex component.
Notice that, by definition\eqref{eq:vfi}, for shocks $s^h$ with non zero
convex compoents, $q(s^h)$ can possibly be zero only when
$h=1$ or $h=l$, i.e. when $s^h$ is the first or the last component of $s$.
In fact,
all other shock components of $s$ are two-sided contact discontinuities
which necessarily must cross at least two connected components of $\MC_k$,
and hence  their strengths are certainly larger
than $2\delta_0$ because of \eqref{eq:delta0}.

Now, defining the {\it inner interaction potential}
of a composed wave $s$ as
\begin{equation}
  \label{eq:iip2}
  Q^I(s) = 2\sum_{h\in\mathcal{S}, \kappa\in\mathcal{R}} \big\vert
  s^h s^\kappa \big\vert + \sum_{h\in\mathcal{S}} q(s^h) +
  \sum_{\kappa\in\mathcal{R}} \big\vert s^\kappa \big\vert^2\,,
\end{equation}
we can finally provide the definition of the
quadratic interaction potential enjoing properties 1-2
by setting
\begin{equation}
  \label{eq:ipd}
  Q_q (t) \doteq 2\sum_{\substack{k_\alpha=k_\beta\\ s_\alpha
  s_\beta>0}} \big\vert s_\alpha s_\beta \big\vert + \,\sum_{\alpha}
  Q^I(s_\alpha) \,+\, c\Bigg[ \sum_{\substack{k_\alpha=k_\beta\\ s_\alpha s_\beta<0}}
  +  \sum_{\substack{k_\alpha<k_\beta\\ x_\alpha(t)>x_\beta(t)}}\Bigg]
  \big\vert s_\alpha s_\beta \big\vert\,,
\end{equation}
where, as usual, $x_\alpha(t)$ denotes the position of the wave $s_\alpha$, and
$k_\alpha$ its characteristic family
while $c$ is the same constant that appearzs in~\eqref{eq:ip}.
Here, the second summation runs
over all composed waves $s_\alpha$ present in~$u^\epsilon(t)$. 
Notice that $Q_q$ differs from the interaction
potential $Q_1$ defined in \S~\ref{sec:oip} only for the presence
of the inner interaction potential $Q^I$ of the
composed waves that replaces the corresponding terms 
of the second and third summands in \eqref{eq:ip1}.
On the other hand, whenever $|s_\alpha|\leq\delta_0$,
we clearly have $Q^I(s_\alpha)=|s_\alpha^h s_\alpha^\kappa|+|s_\alpha^\kappa|^2$,
$h\in\mathcal{S}, \kappa\in\mathcal{R}$,
and thus one recovers the same expression present in $Q_1$.
\begin{remark}
\label{rem:rarcom}
Consider  a shock wave $s$ with strength $|s|\leq\delta_0$ that crosses 
a connected component $\MC_{k,h}$ of $\MC_k$. According
with the above definitions $s$ is decomposed in
a convex and a concave component $s^\smile,\, s^\frown$.
Relying on \cite[Propositions~2.1-2.2]{note}, we deduce that,
choosing $\delta_0$ 
sufficiently small, there holds
\begin{equation}
  \label{eq:srar}
  |s^\smile| \leq c_1\, |s^\frown|\,,
\end{equation}
for some constant $0<c_1<1$.
Such a bound will be useful in the 
study of the variation of the intrinsic
interaction potential $q(s)$ in presence of interactions.
Notice that the above estimate holds even in the case, instead of \eqref{eq:transv-ld},
we assume that there is some even index $p$ so that the following weaker condition
is satisfied:
\begin{equation}
\label{eq:transvp-ld}
D^j_{r_k}\lambda_k(u)=0 \quad \forall~j<p\,,\quad D^p_{r_k}\lambda_k(u)\neq 0\qquad\ \ \forall~u\in\MC_{k,h}\,,
\end{equation}
where $D^j_{r_k}\lambda_k(u)$ denotes the $j$-th derivative of $\lambda_k$ along $r_k$,
inductively defined by setting $D_{r_k}\lambda_k(u)\doteq \nabla\lambda_k(u)\cdot r_k(u)$,
and $D^j_{r_k}\lambda_k(u)\doteq D^{j-1}_{r_k}\lambda_k(u)\cdot r_k(u)$ for all $j>1$.
\end{remark}

Towards an analysis of the interaction
potential above introduced,  we first define a quadratic quantity of
interaction  as in Section~\ref{sec:oip}  for waves of the same family and
with the same sign, to measure the decrease
of the quadratic functional $Q_q$ in \eqref{eq:ipd} 
when waves of this type
with strength  $\leq\delta_0$ are involved in an interaction.
\begin{definition}
\label{def:qid}
Consider two nearby waves
of sizes $s', s''$ with the same sign and belonging to  the the same $k$-th characteristic family, 
with $s'$ located at the left of $s''$.
Assume that $|s'|, |s''|\leq\delta_0$ and, with the same notations of Definition~\ref{def:qi-sld},ù
suppose that the state $u^\sharp$ belongs to the connected component  of $\Omega\setminus\MC_k$
lying between two consecutive manifolds $\MC_{k,h-1}$, $\MC_{k,h}$  of $\MC_k$.
Then, in the case $s', s''>0$,
 we define the {\it quantity of interaction} between $s'$ and $s''$ as
\begin{equation}
  \label{eq:qi2}
  {I} (s', s'')\doteq
  \begin{cases}
    \big[ \big\vert (s' + s'')^r - s'^r \big\vert +
    \vert s'^s \vert \big] \vert s''^s \vert\\
    &\hspace{-2truecm} \text{if}\quad \nabla (\nabla\lambda_k \cdot r_k)(u)\cdot r_k(u)
      \Big\vert_{u\in \MC_{k,h}}< 0\,,\\
    \noalign{\medskip}
    \big[ \big\vert (s' + s'')^r - s''^r \big\vert +
    \vert s''^s \vert \big] \vert s'^s \vert\\
    &\hspace{-2truecm} \text{if}\quad \nabla (\nabla\lambda_k \cdot r_k)(u)\cdot r_k(u)
      \Big\vert_{u\in \MC_{k,h}}> 0\,.
  \end{cases}
\end{equation}
An entirely similar definition is given in the case $s', s''<0$.
For notational convenience we also set $I(s', s'')\doteq 0$
for every pair of waves $s', s''$ of the same family that
have opposite sign.
\end{definition}

Next, following~\cite[Definition~3.5]{sie}, we introduce a definition
of  quantity of interaction 
for a general strictly hyperbolic system \eqref{sys}, which 
measures the decrease of the cubic functional $\Q$ in \eqref{eq:ipngnl} 
when waves of the same family and with the same sign interact together.
\begin{definition}
\label{def:qi}
Consider two nearby waves
of sizes $s', s''$ with the same sign and
 belonging to  the the same $k$-th characteristic family, 
with left states  $u', u''$, respectively.
Let $\widetilde
  F^{\prime} \doteq \widetilde F_k[u'](s',\,\cdot\,)$ and
  $\widetilde F^{\second}\doteq \widetilde
  F_k[u''](s'',\,\cdot\,)$ be the reduced flux with starting point
  $u'$, $u''$, evaluated along the solution of~(\ref{eq2:Ti}) on
  the interval $[0,s^\prime]$, and $[0,s^\second]$, respectively
  (cfr. def.~(\ref{sigmaFdef})).
Then, assuming that~$s'\geq 0$, we say that the {\it amount of interaction} $\J(s^\prime,\,s^{\second})$ between $s^\prime$ and $s^{\second}$ is
the quantity 
     \begin{equation}
      \begin{aligned}
	\label{eq2:qi1}
	\!\!\!&\J(s^\prime,\,s^{\second}) \doteq \int_0^{s^\prime}
	\left\vert \conv_{[0,\,s^\prime]} \widetilde F^{\prime}(\xi) -
	\conv_{[0,\,s^\prime+s^\second]}
	\widetilde F^{\prime}
	\!\cup\! \widetilde F^{\second}
	(\xi) \right\vert d\xi\\
	\noalign{\smallskip}
	&\qquad + \!\!\int_{s^\prime}^{s^\prime+s^\second} \!\!\left\vert
	\widetilde F^{\prime} (s^\prime)+\conv_{[0,\,s^\second]} 
	\widetilde F^{\second} (\xi-s^\prime)
	-\conv_{[0,\,s^\prime+s^\second]} \widetilde F^{\prime}
	\!\cup\! \widetilde F^{\second} (\xi) \right\vert d\xi,
      \end{aligned}
    \end{equation}
    where $\widetilde F^{\prime} \!\cup\! \widetilde
    F^{\second}$ is the function defined on $[0,\,s'+s'']$ as
    \begin{equation}
      \label{eq2:FpFs}
      \widetilde F^{\prime} \!\cup\! \widetilde F^{\second} (s)
      \doteq
      \begin{cases}
	\widetilde F^{\prime}(s) \quad &\text{{\rm if}\quad $s\in
	  [0,s^\prime]$\,,}\\
	\noalign{\smallskip}
	\widetilde F^{\prime}(s^\prime) + \widetilde
	F^{\second}(s-s^\prime) \quad &\text{{\rm if}\quad $s\in
	  [s^\prime, s^\prime+s^\second]$\,.}
      \end{cases}
    \end{equation}
Here, $\conv_{[a,b]} f$, $\conc_{[a,b]} f$ denote the lower convex envelope and the
  upper concave envelope of $f$ on $[a,b]$, defined as in \eqref{convenv}.
  In the case where $s^\prime<0$, one replaces in
  (\ref{eq2:qi1})
  the lower convex envelope with the
  upper concave one, and vice-versa.
As in Definition~\eqref{def:qid}, for notational convenience we also set $\J(s', s'')\doteq 0$
for every pair of waves $s', s''$ of the same family that
have opposite sign.
\end{definition}
\begin{remark}
\label{rem4:amshocks}
Notice that by the Lipschitz continuity of the derivative
$(u,s)\mapsto D_\tau \widetilde F_k[u](s,\cdot)$ 
of the reduced flux \eqref{sigmaFdef}(cfr.~\cite{srp}), it follows 
$\J(s',s'')=\mathcal{O}(1) \cdot |s' s''|$.
Moreover, by Remark~\ref{rem2:consform} one can easily verify that, in the
conservative case, if $s', s''$ are both shocks of the $k$-th family
that have the same sign, then the amount of interaction in
(\ref{eq2:qi1}) takes the form
$$
\J (s', s'')=\big|s' s''\big| \Big|\sigma_k[u^L, u^M]-\sigma_k[u^M,
  u^R]\Big|\,,
$$
i.e. it is precisely the product of the strength of the waves times
the difference of their Rankine Hugoniot speeds.
\end{remark}
Relying on the results in \cite[Section~3]{sie} and on
Lemma~\ref{lem:ie1}, we will  show now that the
interaction potential $Q$ defined by  \eqref{eq:ipngnl}, \eqref{eq:ip}, 
\eqref{eq:ipd}, is 
decreasing at every interaction, and that the variation of the total strength of 
waves $V$ in an approximate solution $u^\eps$
is controlled by $|\Delta Q|$. 
\begin{lemma}
\label{lem:ie}
Under the assumption (H), in the same setting of Lemmas~\ref{lem:ie1}-\ref{lem:ie12} there exists some
constant $c>0$ (in \eqref{eq:ip}, 
\eqref{eq:ipd}), so that
 there
holds
\begin{multline}
  \label{eq:tvest}
  \Delta V  \leq -\sum_{1\leq i\leq N} \C(s'_i,s''_i) +
  \mathcal{O}(1) \cdot \Bigg[ \sum_{\substack{1\leq
  i,j \leq  N\\ \noalign{\smallskip} i> j}} \vert s'_i s''_j \vert 
   + \sum_{1\leq i \leq N}
  \J(s'_i, s''_i) \Bigg]\,,
\end{multline}
\begin{equation}
\begin{aligned}
  \Delta Q \leq &-\frac{1}{2} \Bigg[ \sum_{\substack{1\leq i,j \leq
  N\\ \noalign{\smallskip} i> j}} \vert s'_i s''_j \vert  +
  \sum_{\substack{1\leq i\leq N\\ \noalign{\smallskip} s'_i s''_i<0}} \vert s'_i s''_i \vert
  + \!\!\!\!\sum_{\substack{
   1\leq i\leq N \\ \noalign{\smallskip}
  \vert s'_i\vert, \vert s''_i\vert  \leq
  \delta_0/2}} \!\!\!\! I(s'_i, s''_i)\,+\\
  \noalign{\smallskip}
  &\qquad\qquad +
  \sum_{1\leq i\leq N}\J (s'_i, s''_i)\Bigg] +
  \mathcal{O}(1) \cdot V^-\cdot\!\! \sum_{1\leq i\leq N} \C(s'_i,s''_i)\,.
\end{aligned}
\label{eq:ipest}
\end{equation}
\end{lemma}
\textbf{Proof.}
A proof of the estimate (\ref{eq:tvest}) can be found in \cite{sie}, thus
we will focus our attention on (\ref{eq:ipest}).  For the sake of
simplicity, we shall consider only the case in which the two 
adjacent Riemann problems are solved by a single wave, say $s'$
and $s''$, $s'$ on the left of $s''$. We distinguish three
cases, depending on the strengths of $s'$ and $s''$ and on their
characteristic families.
\begin{enumerate}
  \item
  $s'$ and $s''$ are waves of the $k'$ and $k''< k'$ characteristic families.
  
  \noindent
  To fix the ideas, let $s',s''>0$.  Observe that for every NGNL $k$-family, since condition
  \eqref{eq:transv-ld}
  implies that the characteristic vector field $r_k$ is tranversal to $\MC_k$, by construction
  it follows that the $k$-elementary curves $T_k$ are transversal to each manifold $\MC_{k,h}$.
  As a consequence of this property one can easily verify that, letting $u', u''$ 
  be the left states of $s', s''$, and denoting $u'^+, u''^+$ the left states of 
  the outgoing waves $s_{k'}$ and $s_{k''}$ of
  the $k'$ and $k''$ characteristic families,
  there holds
  \begin{equation}
  \label{eq:tfest1}
  \begin{aligned}
   \Sigma'&\doteq \Bigg[\bigg|\sum_{h\in\SC}
  \sum_{p\in\smile}s'^{h,p}-\sum_{h\in\SC}\sum_{p\in\smile} s_{k'}^{h,p}\bigg| +
  \bigg|\sum_{\kappa\in\RC}s'^\kappa-\sum_{\kappa\in\RC}s_{k'}^\kappa\bigg|\Bigg]
   \\
  &= \mathcal{O}(1) \cdot \Big\Vert T_{k'}[u'] - T_{k'}[u'^+]
  \Big\Vert_{L^\infty}\,,
  \\
  \noalign{\smallskip}
  \Sigma''&\doteq\Bigg[\bigg|\sum_{h\in\SC}
  \sum_{p\in\smile}\!s''^{h,p}-\!\sum_{h\in\SC}\sum_{p\in\smile}  s_{k''}^{h,p}\bigg| +
  \bigg|\sum_{\kappa\in\RC}\!s''^\kappa-\!\sum_{\kappa\in\RC}\!s_{k''}^\kappa\bigg|\Bigg]
  \\
  &= \mathcal{O}(1) \cdot \Big\Vert T_{k''}[u''] - T_{k''}[u''^+]  
  \Big\Vert_{L^\infty}\,.
  \end{aligned}
  \end{equation}
   Here the $L^\infty$ norm in the first and second equality is referred to the intervals
  $[0, \min\{ s', s_{k'} \}]$ and $[0, \min\{ s'', s_{k''} \}]$, respectively. 
  Then, since the interaction estimates in \cite[Section~3]{sie} imply
  \begin{equation}
  \label{eq:rfest1}
   \begin{aligned}
  \Big\Vert T_{k'}[u'^+] - T_{k'}[u']
  \Big\Vert_{L^\infty}
  &= \mathcal{O}(1) \cdot s'',
  \\
  \noalign{\medskip}
  \Big\Vert T_{k''}[u''^+] - T_{k''}[u'']
  \Big\Vert_{L^\infty}
  &=\mathcal{O}(1) \cdot  s'\,,
  \end{aligned}
  \end{equation}
  relying on \eqref{eq:tfest1}-\eqref{eq:rfest1}, 
  we deduce the following bounds on the variation
  of the inner interaction potential at \eqref{eq:iip2}
  \begin{equation}
  \label{eq:innpest1}
  \begin{aligned}
  \Delta Q^I (s')&=  \mathcal{O}(1) \cdot \Sigma'\cdot s'=  \mathcal{O}(1) \cdot s's''\,,
  \\
  \noalign{\medskip} 
  \Delta Q^I (s'')&=  \mathcal{O}(1) \cdot \Sigma'' \cdot s''= \mathcal{O}(1) \cdot s's''\,.
  \end{aligned}
  \end{equation}
  Hence, using
  (\ref{eq:tvest}), \eqref{eq:innpest1}, one obtains
  \begin{equation}
  \label{eq:dQ}
  \Delta Q \leq -c\,\big\vert s' s''\big\vert + \mathcal{O}(1) \cdot
  \vert s's''\vert + \mathcal{O}(1) \cdot \vert s's''\vert \cdot V^-\,,
  \end{equation}
  from which we derive (\ref{eq:ipest}), choosing $c>0$ sufficiently large in \eqref{eq:ipd}.
\smallskip

  \item
  $s'$ and $s''$ are both $k$-waves and $s's''<0$.
  
  \noindent
  By defintion of $Q^I$, and with the same analysis in the previous point,
  one deduces that in this case the inner interaction potential 
  of the outgoing $k$-wave $s$ satisfies 
  $Q^I(s)\leq \min\{Q^I(s'), Q^I(s'')\}+\mathcal{O}(1) \cdot \vert s's''\vert$.
  Hence $\Delta Q^I=\mathcal{O}(1) \cdot \vert s's''\vert$, and thus we obtain 
  the same estimate in~(\ref{eq:dQ}). 
  On the other hand, relying on \cite[Proposition~4.1]{sie} we derive
  \begin{equation}
  \label{eq:dQB1}
  \Delta \Q  \leq\mathcal{O}(1) \cdot V^-\cdot \C(s',s'')\,,
  \end{equation}
  which, together with~(\ref{eq:dQ}),
 yields  (\ref{eq:ipest}), choosing $c>0$ sufficiently large in \eqref{eq:ipd}.
\smallskip

  \item
  $s'$ and $s''$ are both $k$-waves and $s's''>0$.
  
  \noindent
  To fix the ideas, let $s',s''>0$, and call $s$ the  outgoing $k$-wave.
  We shall distinguish a number of cases, depending on 
  the strengths of $s', s''$.
  \begin{enumerate}
  \item
  $\max \{ s', s'' \} \leq \delta_0/2$.
  
  \noindent
  In this case, by definitions \eqref{eq:vfi}-\eqref{eq:ipd} 
  one has $q(s') = q(s'') =$ $=q(s'+s'')=0$, and 
  $\Delta Q_q\leq \Delta q + \Delta Q_1 +\OL (1)\cdot |\Delta V|$, where $\Delta Q_1$
  denotes the variation of the interaction potential $Q_1$ in~\eqref{eq:ip1}
  (related to the waves involved in the interaction).
  Hence, relying on \eqref{eq:tvest}, and  
  applying  (\ref{eq:ip1est}), we deduce $\Delta q=\mathcal{O}(1) \cdot|\Delta V|
  =\mathcal{O}(1) \cdot \J (s',s'')$
  and
  \begin{equation}
   \label{eq:dQq}
  \Delta Q_q  \leq -\frac{1}{2} {I}(s',s'')+\mathcal{O}(1) \cdot \J (s',s'')\,.
  \end{equation}
  On the other hand, due to \cite[Proposition~4.1]{sie}, we get
  \begin{equation}
  \label{eq:dQB}
  \Delta \Q  \leq -\frac{1}{2} \J (s',s'')\,,
  \end{equation}
  which
 yields (\ref{eq:ipest}) choosing $c>0$ 
  sufficiently large in \eqref{eq:ip}.
  \item
  $\delta_0/2<\max \{ s', s'' \} \leq 2 \delta_0$ .
  
  \noindent
  To fix the ideas, assume that $s'$ crosses
   a connected component $\MC_{k,h}$ of $\MC_k$ where
 $\nabla (\nabla\lambda_k \cdot r_k)\cdot r_k< 0$. Because of \eqref{eq:delta0},
this implies that the wave $s''$ on the right of $s'$ must be a shock
with zero convex component and hence $q(s'')=0$.
For sake of simplicity, we shall treat only the case in which also
$s'$ is a shock and $s'\leq \delta_0\leq s/2$, the other cases being similar 
or simpler since for such values of $s', s$ there is the largest possible increase
of $q$ due to the
fact that, by definitions \eqref{eq:vfi}, \eqref{eq:iip1},  one has $q(s')=\varphi(s')=0$,
$\varphi(s)=1$. Under these assumptions, by definitions \eqref{eq:iip1}, \eqref{eq:iip2} we have 
$Q^I(s')=Q^I(s'')=0$, and letting $s^\smile, s^\smile$, $s'^\frown, s'^\frown$,
denote the convex and concave components of $s, s'$,
relying on \eqref{eq:tvest} we deduce
\begin{equation}
\label{eq:qib}
\begin{aligned}
Q^I(s)=q(s)&=2 s^\smile s^\frown + (s^\smile)^2
\\
&=2 s'^\smile (s'^\frown+s'')+(s'^\smile)^2+\OL(1) \cdot |\Delta V|
\\
&\leq 2s'^\smile (s'+s'')+\OL(1) \cdot \J (s',s'')\,.
\end{aligned}
\end{equation}
Moreover, observe that $s'\leq \delta_0\leq s/2$ implies $s'\leq s''+\OL(1) \cdot \J (s',s'')$.
Therefore, using \eqref{eq:qib},
and recalling that by Remark~\ref{rem:rarcom} one has~$s'^\smile<s'^\frown$, we find
\begin{equation}
\label{eq:dQq2}
\begin{aligned}
\Delta Q_q  & \leq-2s's''+Q^I(s)+\OL(1) \cdot |\Delta V|
\\
&\leq-2s''(s'^\frown-s'^\smile)+\OL(1) \cdot \J (s',s'')
\\
&\leq \OL(1) \cdot \J (s',s'')\,.
\end{aligned}
\end{equation}
Hence, \eqref{eq:dQq2} together with \eqref{eq:dQB},
that continues to hold, yields (\ref{eq:ipest}) choosing $c>0$ 
  sufficiently large in \eqref{eq:ip}.

  \item
  $\min \{ s', s'' \}\leq 2 \delta_0 <\max \{ s', s'' \}$.
  
  \noindent
  To fix the ideas assume that $s'$ is a composed wave
  of size $s'\leq 2 \delta_0$, crossing a connected component 
of $\MC_k$ where
 $\nabla (\nabla\lambda_k \cdot r_k)\cdot r_k< 0$. Because of \eqref{eq:delta0},
this implies that the first component $s''^1$ of the (possible composed) wave $s''$ 
on the right of $s'$ must be a shock
of size $s''^1>2 \delta_0$. For sake of simplicity we shall treat only the case in which
also $s'$ is a shock, the other cases being similar.
Observe that, letting $u', u''$ be the left states of $s', s''$,
calling $u'^+$ the left state of the outgoing $k$-wave $s$, 
and letting $s'^\smile, s'^\frown$ denote the convex and concave components
of $s'$,
by the same arguments at point 1 we find
  \begin{equation}
  \label{eq:ttest2}
  \begin{aligned}
   \Sigma&\doteq \Bigg[\bigg|s'^\smile+\sum_{h\in\SC}
  \sum_{p\in\smile}s''^{h,p}-\sum_{h\in\SC}\sum_{p\in\smile} s^{h,p}\bigg| +
  \bigg|\sum_{\kappa\in\RC}s''^\kappa-\sum_{\kappa\in\RC}s^\kappa\bigg|\Bigg]
   \\
  &= \mathcal{O}(1) \cdot \bigg[\Big\Vert T_k[u'] - T_k[u'^+] 
  \Big\Vert_{L^\infty}+\Big\Vert T_k[u''] - T_k[u'^+](s'+\cdot)
  \Big\Vert_{L^\infty}\bigg]\,,
  \end{aligned}
  \end{equation}
where the $L^\infty$ norm of the two terms in the equality is referred to the intervals
$[0, s']$ and  $[0, \min\{ s'', s-s' \}]$, respectively.
Then, applying the interaction estimates in \cite[Section~3]{sie},
we derive
\begin{equation}
\label{eq:rfest3}
\Sigma=\OL(1)\cdot \J(s', s'')\,.
\end{equation}
On the other hand, observe that by definition \eqref{eq:vfi} the above assumptions
imply $\varphi(s''^1)=1$, $\varphi(s^1)=\OL(1)\cdot \J(s', s'')$, since 
the first  component $s^1$ of $s$ satisfies the lower bound
$s^1\geq s''^1+\OL(1)\cdot \J(s', s'')>2\delta_0+\OL(1)\cdot \J(s', s'')$.
Thus, relying on \eqref{eq:tvest}, \eqref{eq:rfest3}, 
and because $s'<s''^1$, we obtain
\begin{equation}
\label{eq:iinpest2}
\begin{aligned}
\Delta Q^I&\leq 2 s'^\smile (s'^\frown + s''^1)+(s'^\smile)^2+\OL(1)\cdot\big(\Sigma+|\Delta V|\big)
\\
&\leq - 2 s''^1(s'^\frown-s'^\smile) + 2 s' s''^1 + \OL(1)\cdot \J(s', s'')\,,
\end{aligned}
\end{equation}
which in turn, recalling that by Remark~\ref{rem:rarcom} one has~$s'^\smile<s'^\frown$, yields
\begin{equation}
\label{eq:dQq3}
\begin{aligned}
\Delta Q_q  & \leq-2s's''^1+\Delta Q^I+\OL(1) \cdot |\Delta V|
\\
&\leq \OL(1) \cdot \big(\J (s',s'')+|\Delta V|\big)
\\
&=\OL(1) \cdot \J (s',s'')\,.
\end{aligned}
\end{equation}
Hence, \eqref{eq:dQq3} together with \eqref{eq:dQB},
that continues to hold, yields (\ref{eq:ipest}) choosing $c>0$ 
  sufficiently large in \eqref{eq:ip}.

  \item
  $\min \{ s', s'' \}> 2\delta_0$.
  
  \noindent
  We shall treat only the case in which the last componenent $s'^{l'}$
  of $s'$ and the first component $s''^1$ of $s''$ are both shocks
  of size $>2\delta_0$. The other cases are simpler or reducible
  to one of the previous cases (a), (b), (c). Then, by definition \eqref{eq:vfi}
  there holds $\vfi (s'^{l'})= \vfi (s''^1)=1$. Moreover, with the
  same notations and with the same arguments of point (c), we have
  \begin{equation}
  \label{eq:ttest3}
  \begin{aligned}
  \Sigma&\doteq 
  \Bigg[\bigg|
  \sum_{h\in\SC}\sum_{p\in\smile}s'^{h,p}+\sum_{h\in\SC}\sum_{p\in\smile}s''^{h,p}-
  \sum_{h\in\SC}\sum_{p\in\smile} s^{h,p}\bigg| +
  \\
  &\qquad\quad +\bigg|
  \sum_{\kappa\in\RC}s'^\kappa+\sum_{\kappa\in\RC}s''^\kappa
  -\sum_{\kappa\in\RC}s^\kappa\bigg|\Bigg]
   \\
  &= \mathcal{O}(1) \cdot \bigg[\Big\Vert T_k[u'] - T_k[u'^+] 
  \Big\Vert_{L^\infty}+\Big\Vert T_k[u''] - T_k[u'^+](s'+\cdot)
  \Big\Vert_{L^\infty}\bigg]\,,
  \\
  &=\OL(1) \cdot \J (s',s'')\,.
  \end{aligned}
  \end{equation}
  Thus, relying on \eqref{eq:tvest}, \eqref{eq:ttest3}, we derive
 \begin{equation}
 \label{eq:iinpest3}
 \begin{aligned}
 \Delta Q_q  & \leq\Delta Q^I+\OL(1) \cdot |\Delta V|
 \\
 &=\OL(1)\cdot\big(\Sigma+|\Delta V|\big)
 \\
 &= \OL(1) \cdot \J (s',s'')\,,
 \end{aligned}
 \end{equation}
  which, together with \eqref{eq:dQB},
that continues to hold, yields (\ref{eq:ipest}) choosing $c>0$ 
  sufficiently large in \eqref{eq:ip}.
  \end{enumerate}
\end{enumerate}
\qed

Relying on the above result one can prove that there exists $C>0$
so that, assuming $V(0)$
sufficiently small, the Glimm functional 
\begin{equation}
  \label{eq:ups}
  t\mapsto \Upsilon (t) \doteq V(t) + C\, Q(t)
\end{equation}
is non increasing at any time, and at every discrete time $t=i\eps$
there holds
\begin{equation}
  \Delta \Upsilon (i\eps) \leq - \frac{1}{2} \Big(\big[ \text{amount of cancellation at
  $t=n\eps$} \big] + \big|\Delta Q(i\eps)\big|\Big)\,.
\end{equation}
Hence, for any given $0\leq m<n$, the total amount of wave interaction and cancellation 
taking place in
the time interval $[m\eps, n\eps]$ is bounded by $\mathcal{O}(1)\cdot|\ups^{m,n}|$,
where
$$
\Delta \ups^{m,n} \doteq \ups (n\eps) - \ups (m\eps)
$$
denotes 
 the variation of $\ups$ on $[m\eps, n\eps]$.

%
%
%	WAVE TRACING
%
%
\section[Wave tracing for general non genuinely nonlinear systems]{Wave tracing for general non genuinely nonlinear systems}
\label{sec:wt}
\setcounter{equation}{0}

We will show now how to implement a wave tracing algorithm
for a NGNL system satisfying the assumption~(H)
so that 
the change in strength and the product of strength 
times the variation in speeds of the primary waves
be bounded by the variation of the Glimm functional in~\eqref{eq:ups}.
Namely, 
recalling the
Definition~\ref{def:par} of a wave partition, we have the following result
analogous to Proposition~\ref{pro:wt1}.
\begin{proposition}
\label{pro:wt2}
Under the assumption (H), the same
conclusions of Proposition~\ref{pro:wt1} hold, with $\Delta\ups^{m,n}$
in place of $\Delta\ups_1^{m,n}$.
\end{proposition}
\textbf{Proof.}
As in the proof of Proposition~\ref{pro:wt1}, in order
to produce a partition for an approximate solution $u^\eps$ that fulfills
properties 1-2,
one may proceed by induction on the time steps $i\eps$, $m\leq i <n$.
Then, assuming that such a partition is given for all times $m\eps \leq t >i\eps$,
our goal is to show how to define a partition of the outgoing waves generated by the 
interactions that take place at $t=i\eps$, preserving the properties 1-2.
As observed in the proof of Proposition~\ref{pro:wt1}, it will
be sufficient to focus our attention on interactions between
waves of the same family and with the same sign, since whenever
any other interaction occurs for a system satisfying the
assumption (H), the change in strength and the product of 
strength times the variation in speeds is controlled by the
variation of a Glimm functional with a quadratic
interaction potential as the part in brackets of \eqref{eq:ipd}
(cfr. \cite[Lemma~3.2 and
Theorem~5.1]{glimmly}).

Thus, consider an interaction between two waves, say $s_k'$, $s_k''$, issuing
from two consecutive mesh points $((i-1)\eps, (j-1)\eps)$ and $((i-1)\eps,
j\eps)$,
 belonging to a $k$-th NGNL characteristic family,
and having the same sign.
Observe that, if $\vert s_k'\vert, \vert s_k''\vert \leq \delta_0/2$, then
relying on the estimates \eqref{eq:tvest}, \eqref{eq:ipest} provided by  Lemma~\ref{lem:ie},
one obtains the desired partition proceeding precisely as in the proof of 
Proposition~\ref{pro:wt1}. Hence, we shall treat only the case 
where $\max \{ \vert s_k'\vert, \vert s_k''\vert \} > \delta_0/2$.
For sake of simplicity, we assume that $s_k',s_k''>0$
and that the outgoing $k$-wave $s_k$ is a shock,
the other cases being entirely similar.
Let
\begin{equation}
\label{partinc2}
\big\{ y'^h_k,\, \lambda'^h_k \big\}_{0<h\leq \ell'}\,,
\qquad\quad 
\big\{ y''^h_k,\, \lambda''^h_k \big\}_{0<h\leq \ell''}\,, 
\end{equation}
be the partitions of $s'_k$ and $s''_k$ enjoing the properties 1-2
(on the interval $[m\eps, (i-1)\eps]$), with sizes
\begin{equation}
  \big\{s'^h_k \big\}_{0<h\leq \ell'}\,,
 \quad\quad 
\big\{s''^h_k\big\}_{0<h\leq \ell''}\,.
\end{equation}
Then, define a partition $\{ y^h_p \}_{0<h\leq \ell_p}$ 
of the outgoing wave
$s_p$ ($p=1,\dots,N$)  
of the $p$-th family 
issuing from $(i\eps, j\eps)$ (with corresponding speeds $\{ \lambda^h_p \}_{0<h\leq \ell_p}$) as in Proposition~\ref{pro:wt1}.
In particular, a partition of~$s_k$ is defined by means of its sizes
as 
\begin{equation}
  \label{eq:par13}
  s^h_k \doteq
  \begin{cases}
    s'^h_k\quad &\text{if\quad
      $h=1,\ldots,\overline\ell'$\,,}\\
    \noalign{\medskip}
    s''^{h-\ell'}_k\quad &\text{if\quad $\overline\ell'=\ell'$\quad and\quad
      $h=\ell'+1,\ldots,\ell'+\overline \ell''$}\\
  \end{cases}
\end{equation}
(with the same notations of the proof of Proposition~\ref{pro:wt1}).

Clearly, such partitions continue to satisfy the bounds \eqref{eq:ttest1},
 \eqref{eq:stest1} and the one-to-one correspondence at \eqref{eq:tcor1},
thanks to the estimate \eqref{eq:tvest}, and because of the inductive 
assumption. Therefore, in order to conclude the proof, it remains to establish
only the estimate \eqref{eq:spest1} on the wave speeds.
To this end, notice that 
the Rankine-Hugoniot speed $\lambda_k$ of the outgoing $k$-wave $s_k$
coincides with
the speeds $\lambda_k^h$ of all subwaves $s^h_k$ defined according with 
Definition~\ref{def:par}, since for
a shock wave the integrand function $\sigma(\cdot)$ in~\eqref{eq:parsp}
results to be a constant (cfr. Remark~\ref{rem2:consform}).
Moreover, by the choice of the speeds of a partition at (\ref{eq:parsp}),
one has
\begin{equation}
\label{eq:par22}
\lambda'^h = \frac{1}{s'^h_k} \int_{\tau'^{h-1}_k}^{\tau'^h_k}
\sigma'(\tau)~d\tau\,, \qquad
\lambda''^h = \frac{1}{s''^h_k} \int_{\tau''^{h-1}_k}^{\tau''^h_k}
\sigma''(\tau)~d\tau\,,
\end{equation}
where $\tau'^h_k \doteq \sum_{p=1}^h s'^p_k$, 
$\tau''^h_k \doteq \sum_{p=1}^h s''^p_k$, and 
$$
\sigma'(\cdot)\doteq \sigma_k[\omega'_k](s'_k,\cdot)\,,\qquad
\sigma''(\cdot)\doteq \sigma_k[\omega''_k](s''_k,\cdot)\,,
$$
denote the map in~\eqref{sigmaFdef} defining the speed of the rarefaction and shock components
of $s'_k$ and $s''_k$, respectively ($\omega'_k, \omega''_k$ being the left
states of $s'_k, s''_k$).

Then, relying on the interaction 
estimates in~\cite[Section~3]{sie}, with the same type of arguments
used in the proof of~\cite[Lemma~3.9]{sie} 
one obtains the following estimate on the wave speeds,
similar to the one provided by
 \cite[Theorem~3.1]{glimmly} :
\begin{equation}
\label{eq:mainwsest}
\begin{aligned}
  {\lambda_k}\cdot \big(s_k'+s''_k\big)&= \int_0^{s'_k+s''_k}
  \sigma(\tau)~d\tau +\mathcal{O}(1) \cdot \J (s'_k,s''_k)\\
  \noalign{\smallskip}
  &=
  \int_0^{s'_k}
  \sigma'(\tau)~d\tau + \int_0^{s''_k}\sigma''(\tau)~d\tau +
  \mathcal{O}(1) \cdot \J (s'_k,s''_k)\\
  \noalign{\smallskip}
  &= 
   \sum_{h=1}^{\ell'}
    s'^h_k \lambda'^h_k + \sum_{h=1}^{\ell''} s''^h_k \lambda''^h_k +
  \mathcal{O}(1) \cdot \J (s'_k,s''_k)\,.
\end{aligned}
\end{equation}
Thus, since by the monotonicity property of $\sigma'(\cdot)$
and $\sigma''(\cdot)$, we have
$$
\lambda''^h_k-\mathcal{O}(1) \cdot
\J (s'_k,s''_k)\leq {\lambda}_k \leq
\lambda'^h_k+\mathcal{O}(1) \cdot
\J (s'_k,s''_k)\qquad\forall~h\,,
$$
using \eqref{eq:mainwsest} we derive
$$
\begin{aligned}
\vert \lambda'^h_k-{\lambda_k}\vert & =
\lambda'^h_k-{\lambda_k} + \mathcal{O}(1) \cdot \J
(s'_k,s''_k)\\
\noalign{\smallskip}
&= \frac{1}{s'_k+s''_k} \cdot \Bigg[ \sum_{\substack{p=1,\ldots,\ell'
}} s'^p_k \big( \lambda'^h_k - \lambda'^p_k \big) + \sum_{p=1}^{\ell''}
  s''^p_k \big( \lambda'^h_k - \lambda''^p_k \big) \Bigg] +\\
\noalign{\smallskip}
&\quad +\mathcal{O}(1) \cdot \frac{\J (s'_k,s''_k)}{s'_k+s''_k}\,,
\end{aligned}
$$
which, in turn, yields
\begin{equation}
\label{eq:stspest1}
\begin{aligned}
  \sum_{h=1}^{\ell'} s'^h_k \vert
  \lambda'^h_k-{\lambda_k}\vert &= \frac{1}{s'_k+s''_k} \cdot \Bigg[
  \sum_{h=1}^{\ell'} \sum_{p=1}^{\ell'} s'^h_k s'^p_k \big( \lambda'^h_k -
  \lambda'^p_k \big) +\\
  \noalign{\smallskip}
  &+\sum_{h=1}^{\ell'} \sum_{p=1}^{\ell''} s'^h_k s''^p_k \big( \lambda'^h_k
  - \lambda''^p_k \big) \Bigg] +\mathcal{O}(1) \cdot \J
  (s'_k,s''_k)\cdot\frac{s'_k}{s'_k+s''_k}.
\end{aligned}
\end{equation}
Notice that 
the terms of the first double sum on the right hand side of~\eqref{eq:stspest1} are
antisimmetric in $(h,p)$, and hence the first summand vanishes. Moreover, recalling
(\ref{eq:par22}), we have
\begin{equation}
\label{eq:stspest2}
\begin{aligned}
  &\sum_{h=1}^{\ell'} \sum_{p=1}^{\ell''} s'^h_k s''^p_k \big( \lambda'^h_k
  - \lambda''^p_k \big) = 
  \\
  &\qquad =\sum_{p=1}^{\ell''} s''^p_k \sum_{h=1}^{\ell'}
  \int_{\tau'^{h-1}_k}^{\tau'^h_k} \sigma'(\tau)~d\tau- \sum_{h=1}^{\ell'} s'^h_k
  \sum_{p=1}^{\ell''} \int_{\tau''^{h-1}_k}^{\tau''^{h-1}_k}
  \sigma''(\eta)~d\eta\allowdisplaybreaks\\
  \noalign{\smallskip}
  &\qquad = s''_k\int_0^{s'_k} \sigma'(\tau)~d\tau - s'_k \int_0^{s''_k}
  \sigma''(\eta)~d\eta\\
  \noalign{\smallskip}
  &\qquad = \int_0^{s'_k}\int_0^{s''_k} \big[ \sigma'(\tau) -
  \sigma''(\eta) \big]~d\eta\, d\tau\,.
\end{aligned}
\end{equation}
On the other hand, observe that the term in \eqref{eq:ipngnl}
corresponding to the outgoing shock wave $s_k$ vanishes
(being the map $\sigma$ constant), and hence one clearly has
\begin{equation}
\label{eq:ipngnld}
\int_0^{s'_k}\int_0^{s''_k} \big| \sigma'(\tau) -
  \sigma''(\eta) \big|~d\eta\, d\tau =
 \mathcal{O}(1) \cdot \big|\Delta \Q (i\eps)\big|\,.
\end{equation}
Thus,
since the assumption $\max \{ s'_k, 
s''_k \}> \delta_0/2$ implies $s'_k+s''_k>\delta_0/2$, 
from \eqref{eq:stspest1}-\eqref{eq:ipngnld} it follows 
\begin{equation}
\label{eq:stspest3}
\begin{aligned}
  \sum_{h=1}^{\ell'}  s'^h_k \vert
  \lambda'^h_k-{\lambda_k}\vert &\leq \frac{2}{\delta_0}
  \int_0^{s'_k}\int_0^{s''_k} \big\vert \sigma'(\tau) - \sigma''(\eta)
  \big\vert~d\eta\, d\tau + \mathcal{O}(1) \cdot \J (s'_k,s''_k)\\
  \noalign{\smallskip}
  &\leq \mathcal{O}(1) \cdot \big|\Delta \Q (i\eps)\big|\,.
\end{aligned}
\end{equation}
An entirely similar estimate can be derived for the components of the
partition of $s''_k$, so that there holds
\begin{equation}
\label{eq:stspest4}
\sum_{h=1}^{\ell''}  s''^h_k\vert 
\lambda''^h_k-{\lambda}_k\vert = \mathcal{O}(1) \cdot \big|\Delta \Q
(i\eps)\big|\,.
\end{equation}
Therefore, relying on the inductive ssumption, from \eqref{eq:stspest3}-\eqref{eq:stspest4}
we recover the desired estimate 
(\ref{eq:spest1}), which completes the proof of the proposition.
\qed

%
%
%	CONCLUSION
%
%
\section[Conclusion]{Conclusion}
\label{sec:con}
\setcounter{equation}{0}

Here we briefly describe how to get the proof of
Theorem~\ref{thm:glimm}, following the ideas contained in \cite{bm}
and relying on the results established in the previous section.

\subsubsection*{Step 1.}
We use the partition of
waves of an approximate solution $u^\eps$ into
$$
\begin{aligned}
  \null
  &\text{primary waves}& &\big\{ \wy^h_k, \wl^h_k \big\}\,,&\\
  \noalign{\smallskip}
  &\text{secondary waves}& &\big\{ \wwy^h_k, \wwl^h_k \big\}\,,&
\end{aligned}
$$
provided by Proposition~\ref{pro:wt2} to construct a piecewise
constant approximation $w=w(t,x)$ of $u^\eps (t,x)$ in a time
interval $[m\eps, n\eps]$ that enjoys the following properties (see
\cite[Section~4]{bm}).
\begin{enumerate}
  \item
    The wave fronts in $w$ are of two kinds, primary and secondary.
  \item
    Each primary front originates at $t=m\eps$ and ends at $t=n\eps$;
  \item
    There is a one-to-one correspondence between primary fronts and
    primary waves $\{ \wy^h_k \}$. In particular, the primary front
    corresponding to $\wy^h_k (m,j)$ has constant size
    $\ws^h_k(m,j)$ and, in view of Proposition~\ref{pro:wt2}, joins
    with a segment the points $(m\eps,j\eps)$ and $(n\eps,
    \ell_{(n,j,k,h)} \eps)$ of the $(t,x)$ plane.
  \item
    The left and right states of the primary front corresponding to
    $\wy^h_k (m,j)$, say $u^{h,L}_k (m,j)$, $u^{h,R}_k (m,j)$, are always
    related by
    $$
    u^{h,R}_k (m,j) = T_k \big[ u^{h,L}_k (m,j) \big] \big( \ws^h_k
    (m,j) \big)\,.
    $$
    Moreover, there holds
    $$
    w(m\eps) = u^\eps (m\eps)\,.
    $$
  \item
    Let $u_\beta^L(t)$, and $u_\beta^R(t)$ be the left and right state of a
    secondary front $x_\beta(t)$ of $w$ at time $t\in [m\eps,n\eps]$. Then, 
    letting {\it CW} denote the set of all pairs of crossing primary waves
    in $u^\eps$ (i.e. all pair of waves $\wy^h_k (m,j),\, \wy^{h'}_{k'} (m,j')$ 
     for which $j<j', k>k'$ and $\ell_{(n,j,k,h)}\geq\ell_{(n,j',k',h')}$),
    there holds
    $$
    \begin{aligned}
    \sum_\beta \big\vert u_\beta^R(t) - u_\beta^L(t) \big\vert &= \mathcal{O}(1) \cdot
    \Bigg[\sum_{j,k,h} \Big\vert \wws^h_k(m,j) \Big\vert +
    \sum_{\textit{CW}}\Big|\widetilde s_k^h(m,j)\,\widetilde s_{k'}^{h'}(m,j')\Big|
    \Bigg]
    \\
    \noalign{\smallskip}
    &= \mathcal{O}(1)
    \cdot \big|\Delta \ups^{m,n}\big|\,,
    \end{aligned}
    $$
    where the summand on the left hand side runs over all secondary fronts in $w(t)$,
    while the second summand on the right hand side runs over all pairs of crossing 
    primary waves in $u^\eps$.
  \item
    All secondary fronts travel with speed $2$, strictly larger than
    all characteristic speeds.
\end{enumerate}

\subsubsection*{Step 2.}
Using the same arguments of \cite[Section~5]{bm}, 
relying on \eqref{eq:disest}, \eqref{smgrbound}, \eqref{eq:stest1}, \eqref{eq:spest1},
one can prove that
\begin{equation}
\label{eq:glimmest}
\begin{aligned}
  &\big\Vert S_{(n-m)\eps} w(m\eps) - w(n\eps)
  \big\Vert_{\elleuno} =\\
  \noalign{\smallskip}
  &\hspace{1in}=\mathcal{O}(1) \cdot \Bigg[ \big|\Delta \ups^{m,n}\big| +
  \frac{1+\log(n-m)}{n-m} + \eps \Bigg] (n-m)\eps\,,\\
  \noalign{\medskip}
  &\qquad \ \ \ \big\Vert u^\eps(n\eps) - w(n\eps) \big\Vert_{\elleuno} =
  \mathcal{O}(1) \cdot \big|\Delta \ups^{m,n}\big| \cdot (n-m)\eps\,,
\end{aligned}
\end{equation}
where $S_{(n-m)\eps} w(m\eps)$ is the semigroup trajectory
related to (\ref{syscon}) with initial datum $w(m\eps)=u^\eps (m\eps)$
evaluated at time $t=(n-m)\eps$.

\subsubsection*{Step 3.}
Now, as in \cite[Section~6]{bm}, let $T=\overline{m} \eps + \eps'$,
for some $\overline m\in\nat$, 
$0\leq \eps'<\eps$, and fix a positive constant $\rho>2\eps$. Them, we
define inductively integers $0=m_0 < m_1 < \ldots < m_\kappa=
\overline{m}$ in this way. Assuming $m_i$ given, then
\begin{enumerate}
  \item
    if $\ups (m_i\eps) - \ups \big( (m_i +1)\eps \big) \leq \rho$, 
    let $m_{i+1}$ be the largest integer less or equal to
    $\overline{m}$ such that $(m_{i+1} - m_i)\eps \leq \rho$ and $\ups
    (m_i\eps) - \ups ( m_{i +1}\eps ) \leq \rho$;
  \item
    if $\ups (m_i\eps) - \ups \big( (m_i +1)\eps \big) >\rho$,  set
    $m_{i+1} \doteq m_i+1$.
\end{enumerate}
On every interval $[m_i\eps,m_{i+1}\eps]$ where 1. holds, we construct a piecewise constant
approximation of $u^\eps$ according to Step~1, and using \eqref{eq:glimmest}
we derive
\begin{equation}
\label{eq:finest1}
\begin{aligned}
&\big\Vert u^\eps(m_{i+1}\eps) - S_{(m_{i+1}-m_i)\eps} u^\eps(m_i\eps) \big\Vert_{\elleuno}
\\
&\qquad\quad= \mathcal{O}(1) \cdot \Bigg[ \big|\Delta \ups^{m_i,m_i+1}\big| +
  \frac{1+\log(m_{i+1}-m_i)}{m_{i+1}-m_i} + \eps \Bigg] \big(m_{i+1}-mi\big)\eps\,.
\end{aligned}
\end{equation}
On the other hand, on each interval $[m_i\eps,m_{i+1}\eps]$ where 2. is verified,
by the Lipschitz continuity of $u^\eps$ and applying \eqref{smgrbound}
we find
\begin{equation}
\label{eq:finest2}
\big\Vert u^\eps(m_{i+1}\eps) - S_{(m_{i+1}-m_i)\eps} u^\eps(m_i\eps) \big\Vert_{\elleuno}
= \mathcal{O}(1) \cdot \eps\,.
\end{equation}
Hence, observing that the cardinality of both classes of intervals 1.-2.
is bounded by 
$\mathcal{O}(1) \cdot \rho^{-1}$, 
from \eqref{eq:finest1}-\eqref{eq:finest2} we finally deduce
$$
\big\Vert u^\eps (T) - S_T \overline{u} \big\Vert_{\elleuno} =
\mathcal{O}(1) \cdot \Bigg[ \rho + \frac{\eps}{\rho} \log
  \frac{\rho}{\eps} + \eps \left( 1+\frac{1}{\rho} \right) \Bigg]\,,
$$
which yields (\ref{eq:Glimm}) choosing $\rho \doteq
\sqrt\eps\cdot \log \vert \log \eps \vert$.
\qed
\begin{remark}
\label{rem8:genngnl}
By the same observations of Remark~\ref{rem3:c11ngnl},
one deduces that
the conclusion of Theorem~\ref{thm:glimm} remains valid 
if we assume that the flux function $F$ is $\C^{2,1}$
and that, for each $k$-th NGNL characteristic family,
the linearly degenerate manifold $\MC_k$ in \eqref{LDM} 
is the union of
a finite number of 
connected manifolds $\mathcal{M}_{k,h}$, 
that are either $N\!-\!1$-dimensional as in (H),
or $N$-dimensional and in this case the following conditions hold.
The flux function  $F$ is~$\C^3$ on $\Omega\setminus\MC_k$,
the vector field $r_k$ is transversal to the boundary of $\MC_{k,h}$,
and letting $\partial^+ \MC_{k,h},\, \partial^- \MC_{k,h}$ denote the connected components of
the boundary of $\MC_{k,h}$
where $r_k$ points towards $\Omega\setminus\MC_{k,h}$ and~$\MC_{k,h}$, respectively,
the one-sided second derivatives of $\lambda_k$ 
(see Remark~\ref{rem3:c11ngnl}) satisfy
\begin{equation}
\begin{aligned}
   &\nabla^+ (\nabla\lambda_k \cdot r_k)(u)\cdot r_k(u)< 0 
\qquad \quad \forall u\in \partial^+\mathcal{M}_{k}\,,
\\
   &\nabla^- (\nabla\lambda_k \cdot r_k)(u)\cdot r_k(u)< 0 
\qquad \quad \forall u\in \partial^-\mathcal{M}_{k}\,,
\end{aligned}
\end{equation}
or the opposite inequalities.
Indeed, if we again add the term $2\,\sum_{\substack{p,q\\ p\neq q}} |s_{\alpha,p}^s\,s_{\alpha,q}^s|$
in the interaction potential $Q_q$ defined by~\eqref{eq:iip2}-\eqref{eq:ipd},
for every wave $s_\alpha$ containing
several contact discontinuities $s_{\alpha,1}^s, \dots ,  s_{\alpha,l}^s$,
the estimates stated in Lemma~\ref{lem:ie} continue to hold,
provided that
\begin{equation}
\inf\Big\{s>0\,:\,R_k[u](s)\in\partial^+ \MC_{k,h}\,,\ u\in\partial^- \MC_{k,h}\Big\}>0\,,
\end{equation}
which is certainly true up to a possible slight restriction of the domain $\Omega$.
Relying on Lemma~\ref{lem:ie}, one then  establishes  Proposition~\ref{pro:wt2} 
with the same arguments of Section~\ref{sec:wt}, and thus conclude as above.
\end{remark}

\end{document}